%% file: sympisotopy.tex
\documentclass{amsart}
\usepackage{a4}
 \usepackage{amsfonts,amsmath,amscd,amssymb,amssymb,amsbsy,amsthm,amstext,amsopn}
 \usepackage{bm,mathrsfs,verbatim,multicol}
\usepackage[dvips]{graphicx,color}
\usepackage{latexsym}
\usepackage{euscript,epsfig}
\usepackage[arrow,matrix,line]{xy}

 \numberwithin{equation}{section}
 \newtheorem{theorem}{Theorem}[section]
 \newtheorem{proposition}[theorem]{Proposition}
 \newtheorem{lemma}[theorem]{Lemma}
 \newtheorem{corollary}[theorem]{Corollary}

 \theoremstyle{definition}
 \newtheorem{definition}[theorem]{Definition}
 
 \newtheorem{remark}[theorem]{Remark}

\input xy
\xyoption{all}

 \newcommand{\CC}{\ensuremath{\mathbb{C}}}
 
 \newcommand{\NN}{\ensuremath{\mathbb{N}}}
 
 \newcommand{\RR}{\ensuremath{\mathbb{R}}}
 \newcommand{\FF}{\ensuremath{\mathbb{F}}}
 \newcommand{\CP}{\ensuremath{\mathbb{CP}}}
 \newcommand{\XX}{\ensuremath{\mathcal{X}}}
 
 \newcommand{\VV}{\ensuremath{\mathcal{V}}}
 \newcommand{\SF}{\ensuremath{\mathbb{S}}}
\newcommand{\crit}{\mathrm{crit}}
\newcommand{\std}{\mathrm{std}}
\newcommand{\loc}{\mathrm{loc}}

\newcommand{\git}{\ensuremath{/\!\!/}}
\newcommand{\pr}{\mathrm{pr}}
\newcommand{\id}{\mathrm{id}}
\newcommand{\s}{\underline{s}}
\newcommand{\x}{\underline{x}}

 \DeclareMathOperator{\Span}{span}

\DeclareMathOperator{\Symp}{Symp}
\DeclareMathOperator{\Diff}{Diff}

\begin{document}

\title{A Symplectic Isotopy of a Dehn Twist on $\CP^n\times\CP^{n+1}$}
\author{Emiko Dupont}
\address{Department of Mathematics, Stony Brook
   University, Stony Brook, NY 11794, USA}

\bibliographystyle{alpha}

\maketitle

\begin{abstract}
The complex manifold $\CP^n\times\CP^{n+1}$ with symplectic form
$\sigma_\mu=\sigma_{\CP^n}+\mu\sigma_{\CP^{n+1}}$, where
$\sigma_{\CP^n}$ and $\sigma_{\CP^{n+1}}$ are normalized Fubini-Study
forms, $n\in\NN$ and
$\mu>1$ a real number, contains a natural Lagrangian sphere
$L^{\mu}$. We prove that the Dehn twist along $L^{\mu}$ is
symplectically isotopic to
the identity for all $\mu>1$. This isotopy can be chosen so that it
pointwise fixes a complex hypersurface in  $\CP^n\times\CP^{n+1}$ and
lifts to the blow-up of $\CP^n\times\CP^{n+1}$ along a complex $n$-dimensional
  submanifold.
\end{abstract}


\section{Introduction}\label{chap:intro}
Suppose $L$ is a Lagrangian sphere in a symplectic manifold
$(M,\omega)$. The generalized Dehn twist (the Dehn
twist for short) $\tau_L$ along $L$ is a symplectomorphism on
$M$ that is compactly supported near $L$ and restricts to $L$ as the
antipodal map. A symplectomorphism $\tau$ is
said to be symplectically isotopic to the identity if there is a
smooth family $(\tau_{s})_{0\le s\le 1}$ of symplectomorphisms
such that $\tau_0=\id$ and $\tau_1=\tau$.

In this paper we consider the (real) $4n+2$-dimensional manifold
$\CP^n\times\CP^{n+1}$ with the product symplectic form
$\sigma_{\mu}=\sigma_{\CP^n}+\mu\sigma_{\CP^{n+1}}$ for $\mu>1$ a real
number. Here $n\in\NN$ and $\sigma_{\CP^n}$ and $\sigma_{\CP^{n+1}}$
are the Fubini-Study forms on $\CP^n$ and $\CP^{n+1}$, respectively,
normalized to integrate to $\pi$ on $\CP^1$. As described in Lemma \ref{lem:hopf}, the
graph of the complex conjugate of the Hopf map $S^{2n+1}\rightarrow\CP^n$
embeds as the Lagrangian sphere 
\begin{equation}\label{Lmu}
L^{\mu}=\{([\overline{z}_0:\cdots:\overline{z}_n],[z_0:\cdots:z_n:\sqrt{\mu-1}])\mid\vert
z\vert^2=1\}
\end{equation} 
in $(\CP^n\times\CP^{n+1},\sigma_{\mu})$. The main result of this paper
is the following.
\begin{theorem}\label{thm:isotopy}
For all $\mu>1$, the Dehn twist $\tau_{L^{\mu}}\in
\Symp(\CP^n\times\CP^{n+1},\sigma_{\mu})$ along the
Lagrangian $L^{\mu}$ is symplectically isotopic to the identity by
an isotopy whose restriction to the complex hypersurface
\[
S=\{([s_0:\cdots:s_n],[x_0:\cdots:x_{n+1}])\in\CP^n\times\CP^{n+1}\mid
s_0x_0+\cdots+s_nx_n=0\}
\]
is the identity.
\end{theorem}
The isotopy of Theorem \ref{thm:isotopy} was established in the
case $n=1$ and $\mu\gg 1$ by Corti-Smith \cite[Section
7]{CS05}. Corti-Smith construct a singular fibration with
non-singular fibers isotopic to $(\CP^1\times\CP^2,\sigma_{\mu})$ and
such that the monodromy around the only singular fiber of this fibration is symplectically isotopic to the Dehn
twist $\tau_{L^{\mu}}$. As the monodromy is also known to be symplectically trivial,
this proves the result. The proof of Theorem \ref{thm:isotopy} uses
exactly the same fibration but by examining the construction in more
detail, we are able to establish an isotopy for all
$\mu>1$. Furthermore, we observe that the construction generalizes to
all $n\in\NN$ and that the complex hypersurface $S$ is pointwise-fixed under
the isotopy. 

The hypersurface $S$ is a $\CP^n$-bundle over
$\CP^n$ in which the base coordinates are $[s_0:\cdots:s_n]$, and each fiber is a
linearly embedded copy of $\CP^n$ in $\CP^{n+1}$. Write $\CP^{n+1}$
as $\CC^n\sqcup D$ where $\CC^n$ is the coordinate chart centered
at the point $p_0=[0:\cdots:0:1]$ and $D=(x_{n+1}=0)\cong\CP^n$. Then $S$ contains a section at $0\in\CC^n$ given by
$S_0=\CP^n\times\{p_0\}$ and a section at infinity,
namely,
\[S_{\infty}=\left\{([s_0:\cdots:s_n],[x_0:x_1:0:\cdots:0])\mid
s_0x_0+s_1x_1=0\right\}\subset \CP^n\times D.
\]
We show that $S_0$ and $S_{\infty}$ are not only pointwise-fixed
under the isotopy but in fact we have the following result.
\begin{corollary}\label{cor:blow-up}
Let $S_0$ and $S_{\infty}$ be the complex submanifolds in $S$, each
isomorphic to $\CP^n$,  defined
by
\[
\begin{array}{ccl}S_0&=&\CP^n\times\{p_0\},\\
S_{\infty}&=&\{([s_0:\cdots:s_n],[x_0:x_1:0:\cdots:0])\in\CP^n\times\CP^{n+1}\mid
s_0x_0+s_1x_1=0\}
\end{array}
\]
where $p_0=[0:\cdots:0:1]$.
The isotopy of Theorem \ref{thm:isotopy} lifts to the blow-up of
$\CP^n\times\CP^{n+1}$ along $S_0$
and $S_{\infty}$ if the size of the
blow-up is sufficiently small.
\end{corollary}
One should note that the Dehn twist $\tau_{L^{\mu}}$ is well-defined
on the blow-up since both $S_0$ and $S_{\infty}$ are disjoint from
$L^{\mu}$ and the Dehn twist is supported near $L^{\mu}$.

\medskip

Theorem \ref{thm:isotopy} is a consequence of
Theorem \ref{thm:lef}. The latter constructs a Lefschetz fibration
with exactly one critical point whose non-singular fibers are isotopic
to $(\CP^n\times\CP^{n+1},\sigma_{\mu})$ and with vanishing cycle
isotopic to the Lagrangian sphere $L^{\mu}$. Thus, the
monodromy around a positively-oriented loop in $\CP^1$ that circles
the critical value exactly once, is symplectically isotopic to the
Dehn twist $\tau_{L^{\mu}}$. Since there is only one critical value, the
monodromy is also symplectically isotopic to the identity.
\begin{theorem}\label{thm:lef}
Let $\mu>1$. There is a Lefschetz fibration $(\XX,\pi,
J,\Omega^{\mu})$ such that:
\begin{enumerate}
\item $\pi$ has exactly one critical point $z_{\crit}\in \XX_0=\pi^{-1}(0)$;\label{lefdelta1}
\item\label{lefdelta2}there is a holomorphic trivialization
\[
\Phi:\XX\setminus\XX_{0}\rightarrow(\CP^n\times\CP^{n+1})\times\CC
\]
whose restriction to
  the fiber $\XX_{\infty}$ over the point at infinity in $\CP^1$
  composed with the projection to $\CP^n\times\CP^{n+1}$ is a
  symplectomorphism
\[
\Phi_{\infty}:(\XX_{\infty},\Omega^{\mu}\mid_{\XX_{\infty}})\rightarrow(\CP^n\times\CP^{n+1},\sigma_{\mu});
\]
\item the map $\Phi_{\infty}$ takes the vanishing cycle
  $\VV_{\infty}^{\mu}$ in $\XX_{\infty}$ to
  $L^{\mu}$ in $\CP^n\times\CP^{n+1}$.\label{lefdelta3}
\end{enumerate}
\end{theorem}
For $n=1$, Theorem \ref{thm:lef} is an extended version of
\cite[Proposition 1]{CS05}. The Lefschetz fibration
$(\XX,\pi,J,\Omega^{\mu})$ is an extension of the
fibration constructed in \cite{CS05} from a disk $\Delta\subset\CC$ to
$\CP^1=\CC\cup\{\infty\}$. Note, however, that unless
$\mu\gg 1$, the full statement of Theorem \ref{thm:lef} is necessary
in order to identify the vanishing cycle with the Lagrangian
$3$-sphere $L^{\mu}$ in $\CP^1\times\CP^2$.  After deforming $\Omega^{\mu}$ to a form
$(\Omega^{\mu})'$ which is standard in a small neighborhood of the
critical point $z_{\crit}$, Corti-Smith compute the vanishing cycle
$(\VV_t^{\mu})'$ in the fibers close to the singular fiber $\XX_0$.
For each $t\in\Delta$, write
$\Phi_t:\XX_t\rightarrow\CP^1\times\CP^2$ for the restriction
$\Phi\mid_{\XX_t}$ of the map $\Phi$ described in
Theorem~\ref{thm:lef} followed by projection to $\CP^1\times\CP^2$.
The map $\Phi_t$ is biholomorphic and
$\alpha_t'=(\Phi_t^{-1})^*((\Omega^{\mu})'\mid_{\XX_t})$ is a
symplectic form on $\CP^1\times\CP^2$ in the same cohomology class
as $\sigma_{\mu}$. Hence by Moser's Theorem, $\alpha_t'$ is isotopic
to $\sigma_{\mu}$. If $\mu\gg 1$, the fiber $\XX_1$ intersects the
neighborhood in which $(\Omega^{\mu})'$ is standard and Corti-Smith
show that the pull-back by $\Phi_1$ of the vanishing cycle $(\VV_1^{\mu})'$ in
$\XX_1$ equals the Lagrangian $L^{\mu}\cong S^3$. Hence
the isotopy from $\alpha_1'$ to $\sigma_{\mu}$ can be made to fix
$L^{\mu}$. 

For general $\mu>1$, however, it is not clear what
happens to the vanishing cycle under the isotopy and, in particular,
we cannot identify it with $L^{\mu}$. On the other hand, by
extending the fibration over $\CP^1$ as in Theorem \ref{thm:lef}, we see that the
form $\alpha_{\infty}=(\Phi_{\infty}^{-1})^*(\Omega^{\mu})$ is not
only isotopic to $\sigma_{\mu}$ but, in fact,
$\alpha_{\infty}=\sigma_{\mu}$ and
$\Phi_{\infty}:\XX_{\infty}\rightarrow\CP^1\times\CP^2$ is a
symplectomorphism. The key step in our construction is the following. By symplectically embedding a large
part of the total space $(\XX, \Omega^{\mu})$ of the Lefschetz
fibration into a toric manifold $(\FF,\omega_{\FF}^{(1,\mu)})$, we
obtain a Darboux chart on $(\XX,\Omega^{\mu})$ that enables us to
compute the horizontal spaces of the symplectic connection coming
from $\Omega^{\mu}$ in a large neighborhood of the critical point
$z_{\crit}$. In particular, we are able to compute the vanishing
cycle $\VV_{\infty}^{\mu}$ explicitly and to see that
$\Phi_{\infty}(\VV_{\infty}^{\mu})$ is the Lagrangian $L^{\mu}$.

Our construction generalizes easily from the case $n=1$ to arbitrary
$n\in\NN$.  Furthermore, with $S\subset\CP^n\times\CP^{n+1}$ as defined
in Theorem \ref{thm:isotopy}, we see that the corresponding
hypersurface $\Phi_t^{-1}(S)$ in the fiber $\XX_t$ is held fixed under
symplectic parallel transport in
$(\XX,\pi,J,\Omega^{\mu})$. This implies that the
isotopy of Theorem \ref{thm:isotopy} fixes the complex hypersurface $S$ in $\CP^n\times\CP^{n+1}$.

\medskip

The results in this paper are motivated by the symplectic isotopy
problem which is described below.
Given a symplectic manifold $(M,\omega)$, the group of
symplectomorphisms on $M$ is an infinite-dimensional Lie
subgroup of the group $\Diff^+(M)$ of orientation-preserving
diffeomorphisms of $M$. We say that a symplectomorphism $\tau$ is
smoothly isotopic to the identity if there is an isotopy
$(\tau_s)_{0\le s\le 1}$ with 
$\tau_s\in \Diff^+(M)$. A symplectomorphism is called essential if it
is smoothly isotopic to the identity but not symplectically so.

\noindent{\bf The symplectic isotopy problem}:
\emph{ Given a symplectic manifold $(M,\omega)$, does it admit
  essential symplectomorphisms?\footnote{In this paper we consider only compact symplectic manifolds.}}

When the dimension of $M$ is $2$, the
space of symplectic forms on $M$ is convex and it follows from
Moser's Theorem that no
essential symplectomorphisms exist. In dimension $4$, Gromov and Abreu-McDuff  (see
\cite[p.320-321]{MS04}) show that the answer to the isotopy problem is
once again negative
for $\CP^2$, $\CP^1\times\CP^1$ and the one point blow-up of
$\CP^2$. In 1997, however, Seidel \cite{Sei97} shows that under fairly
weak conditions on a $4$-manifold, if $\tau_L$ is the Dehn twist in a
Lagrangian sphere $L$, then $\tau_L^2$ is an essential
symplectomorphism. Notice here that since $\tau_L$ restricts to $L$ as
the antipodal map, $\tau_L$ cannot be isotopic to the
identity when the dimension of $L$ is even. On the other hand, Seidel
showed that when $M$ is $4$-dimensional, the square of the Dehn twist
is always smoothly isotopic to the identity.

In dimensions $6$ and above, very little is known. It is not even clear
that Dehn
twists (or their squares) are smoothly isotopic to the
identity. Notice that in the
example $(\CP^n\times\CP^{n+1},\sigma_{\mu})$ studied in this paper, the Lagrangian sphere is
odd-dimensional. Thus, there is no homological obstruction for the
Dehn twist to be isotopic to the identity, and indeed, Theorem
\ref{thm:isotopy} shows that it is. Although Theorem \ref{thm:isotopy} does not provide an answer to the symplectic
isotopy problem, the construction of the isotopy relies on certain
symmetries specific to these examples that suggests they may be
higher-dimensional analogues of the special cases $\CP^2$ and
$\CP^1\times \CP^1$ in dimension $4$. (See \cite[Concluding remarks]{ED}). 

Furthermore, by
understanding how the isotopy could be destroyed, one may be able to
construct examples of essential symplectomorphisms in these dimensions.
In contrast to Corollary \ref{cor:blow-up}, the proof of \cite[Proposition 2]{CS05} shows the
following result about $S_0$ and a different $n$-dimensional
submanifold $S_{\infty}'$ at infinity.
\begin{proposition}\label{prop:cannotfix}
Consider the complex submanifolds in $\CP^n\times\CP^{n+1}$ given by
\[\begin{array}{ccc}
S_0=\CP^n\times\{p_0\}&\textrm{and}&
S_{\infty}'=\CP^n\times\{p_{\infty}\},
\end{array}
\]
where $p_0=[0:\cdots:0:1]$ and $p_{\infty}=[1:0:\cdots:0]$.
There is no smooth isotopy
between the Dehn twist $\tau_{L^{\mu}}$ and the identity that
simultaneously fixes $S_0$ and $S_{\infty}'$ pointwise.
\end{proposition}
It therefore seems likely that blowing up along the
submanifolds $S_0$ and $S_{\infty}'$ would destroy the symplectic
isotopy of Theorem \ref{thm:isotopy}. This does not preclude there
being no symplectic isotopies between the Dehn twist along $L^{\mu}$
and the identity in the blow-up. If one could nevertheless
exhibit a smooth isotopy between $\tau_{L^{\mu}}$ and the identity in
the blow-up, then this example would provide a candidate for an
essential symplectomorphism in dimension $4n+2$. An alternative approach to
destroying the isotopy of Theorem
\ref{thm:isotopy} would be to perform a large blow-up of the
submanifolds $S_0$ and $S_{\infty}$ of Corollary \ref{cor:blow-up}. It would
be interesting to estimate the bounds of the size of
the blow-ups allowed for the isotopy to persist and to see what
happens when we exceed the bounds.


\subsection*{Acknowledgments}
This paper contains the results of my Ph.D. thesis, and I
thank my Ph.D. advisor, Dusa McDuff, for her help and guidance
throughout my Ph.D. I am also grateful to Aleksey Zinger for his
many comments and suggestions.


\section{Preliminaries}\label{chap:lef-dehn}
This section is expository and can be skipped if one is already
familiar with Dehn twists and Lefschetz fibrations. The main statement is Proposition
\ref{prop:lef-dehn} which shows that the monodromy of a Lefschetz
fibration around a loop that circles a critical value once is
symplectically isotopic to the Dehn twist along the vanishing cycle
corresponding to the critical value. The discussion in this section is based on  \cite{Sei97} and
\cite[Section 1]{Sei03}. Note that although \cite{Sei03} assumes exactness of
Lefschetz fibrations, i.e., that each non-singular fiber is a
symplectic manifold with boundary and the symplectic form on these fibers is exact,
the proofs of the results that we use are easily adapted to ordinary
Lefschetz fibrations. Details can be found in \cite{ED}.

\subsection{The model Dehn twist}\label{sec:model}
Consider the cotangent bundle of $S^N$
\[
T^*S^N=\{(u,v)\in\RR^{N+1}\times\RR^{N+1}\mid\Vert u\Vert=1,\;\langle u,v\rangle=0\}
\]
with symplectic form $\omega_{T^*S^N}=-\sum_j du_j\wedge
dv_j$. The zero section $L_0$ is a Lagrangian submanifold of
$(T^*S^N,\omega_{T^*S^N})$. The length function $h:T^*S^N\rightarrow
\RR$ given by $h(u,v)=\Vert v\Vert$ generates a Hamiltonian circle
action on $T^*S^N\setminus L_0$ whose flow is given by
\begin{equation}
\varphi^h_{\theta}(u,v)=\left(\cos(\theta)v-\sin(\theta)\Vert v\Vert
u,
  \cos(\theta)u+\sin(\theta)\tfrac{v}{\Vert v\Vert}\right).
\end{equation}
The time-$\pi$ map $\varphi^h_{\pi}$ extends over the zero section
by the antipodal map $A(u,0)=(-u,0)$. Now let $R:\RR\rightarrow\RR$
be a smooth function such that $R(s)=0$ for $s\ge s_R$ for some
$s_R>0$ and $R(-s)=R(s)-s$ for all $s$. Let $H=R\circ h$. The flow of $H$ is
$\varphi^H_{\theta}(u,v)=\varphi^h_{\theta R'(\Vert v\Vert)}(u,v)$. Since $R'(0)=\frac{1}{2}$, $\varphi^H_{2\pi}$ extends continuously
to $T^*S^N$ by the antipodal map. By \cite[Lemma 1.8]{Sei03} this
extension is smooth and hence is a symplectomorphism.
\begin{definition}\label{def:modeldehntwist}
Let $\tau$ be the time-$2\pi$ flow of $H=R\circ h$ on $T^*S^N\setminus
L_0$ as above, extended to $L_0$ by the antipodal map. The symplectomorphism
$\tau$ is called a \emph{model Dehn twist}.
\end{definition}

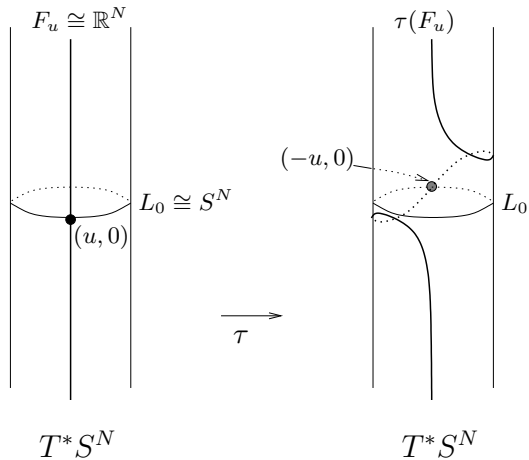
\begin{figure}[!ht]
\begin{center}
\input{dehntwistuv.pstex_t}
\caption{The model Dehn twist\label{fig:dehntwist}}
\end{center}
\end{figure}

Figure \ref{fig:dehntwist} shows the image of a
fiber $F_u=\pi_{T^*S^N}^{-1}(u,0)$ of $T^*S^N$ under a model Dehn twist. Here
$\pi_{T^*S^N}:T^*S^N\rightarrow L_0\cong S^N$ is the natural
projection. We see that $(u,v)\in T^*S^N$ with $\Vert v\Vert$
large are held fixed by $\tau$, but as $\Vert v\Vert$ decreases, the
circle action $\varphi^h_{\theta}$ is ``turned on''
 with larger and larger $\theta\in (\pi,0)$ until at $\Vert
 v\Vert=0$, we reach the antipodal map $\varphi^h_{\pi}$. Although the definition of the model Dehn twist
depends on the choice of function $R$, the symplectic isotopy class
of $\tau$ is independent of $R$.

Let $\omega_{\CC^N}=\frac{i}{2}\sum_jdz_j\wedge d\overline{z}_j$
denote the standard symplectic structure on $\CC^N$. Consider the singular
fibration $\pi_{\std}:\CC^N\rightarrow\CC$ given by
\[
\pi_{\std}(z_1,\ldots,z_N)=\sum_{j=1}^N z_j^2.
\]
All fibers of $\pi_{\std}$, except for
$\pi_{\std}^{-1}(0)$, are symplectic submanifolds of
$(\CC^N,\omega_{\CC^N})$. For each $r>0$, 
\[
\pi_{\std}^{-1}(r):=\{z\in\CC^N\mid\Vert \Re (z)\Vert^2-\Vert \Im
(z)\Vert^2=r,\;\sum_j\Re (z_j)\Im (z_j)=0\}
\]
is symplectomorphic to $(T^*S^{N-1},\omega_{T^*S^{N-1}})$ by the map
\begin{equation}\label{Phistd_r}
\Phi_r(z)=\left(\frac{\Re (z)}{\Vert \Re
(z)\Vert},-\Vert \Re (z)\Vert \Im (z)\right).
\end{equation}

Away from the singular fiber, the symplectic structure $\omega_{\CC^N}$ gives a natural connection
on $\pi_{\std}:\CC^N\rightarrow\CC$ for which parallel transport maps
are symplectomorphisms. Each non-singular fiber $\pi_{\std}^{-1}(t)$, $t\ne 0$, contains a Lagrangian sphere
\begin{equation}\label{Sigma_t}
\Sigma_t=\{\sqrt{r}e^{i\frac{\theta}{2}}z\in\CC^N\mid \Im
(z)=0,\;\Vert z\Vert^2=1\},\quad \textrm{where }t=re^{i\theta}.
\end{equation}
These spheres are mapped to each other by parallel transport. Note that for $r>0$, $\Phi_r(\Sigma_r)$ is the zero section
$L_0$. 
The following result
is based on the proof of \cite[Lemma 1.10]{Sei03} and is proved in detail
in \cite{ED}.
\begin{lemma}\label{lem:stdmodel}
Consider the singular symplectic fibration
$\pi_{\std}:\CC^N\rightarrow\CC$ given by
$\pi_{\std}(z_1,\ldots,z_N)=\sum_{j=1}^N z_j^2$, where $\CC^N$ has
the standard symplectic form $\omega_{\CC^N}$. For $r>0$, let
$\rho_r:\pi^{-1}_{\std}(r)\rightarrow\pi^{-1}_{\std}(r)$ denote the
monodromy around the loop $\gamma(\theta)=re^{i\theta}$,
$0\le\theta\le 2\pi$, and let
$\Phi_r:\pi^{-1}_{\std}(r)\rightarrow T^*S^{N-1}$ be defined
as in (\ref{Phistd_r}). Then for all $r>0$,
\[
\widetilde{\tau_r}=\Phi_r\circ\rho_r\circ(\Phi_r)^{-1}:
T^*S^{N-1}\rightarrow T^*S^{N-1}
\]
is symplectically isotopic to a model Dehn twist.
\end{lemma}

\subsection{Dehn twists}
\begin{definition}\label{def:dehntwist}
Let $(M^{2N},\omega)$ be a symplectic manifold and $L\subset M$ a
Lagrangian sphere with a chosen identification $\iota:S^N\rightarrow
L$. By the Lagrangian Neighborhood Theorem, for some small
$\lambda>0$, we can extend $\iota$ to a symplectic embedding of the space $T^*_{\le \lambda}S^N$ of
cotangent vectors with $\Vert v\Vert \le\lambda$ to a neighborhood
$\mathcal{N}(L)$ of $L$. Choose a model Dehn twist $\tau$
whose support is in $T^*_{\le\frac{\lambda}{2}}S^N$. The
\emph{generalized Dehn twist} (or simply \emph{Dehn twist}) along $L$ is the map
\[
\tau_L(z)=\begin{cases}\iota\circ\tau\circ\iota^{-1},&\textrm{if
    }z\in\mathcal{N}(L);\\
z&\textrm{otherwise.}
\end{cases}
\]
\end{definition}

\begin{remark}\label{framing}
The Dehn twist is a symplectomorphism of $M$ that
depends on the choice of identification $\iota:S^N\rightarrow L$. We
say that two identifications $\iota_1$ and $\iota_2$ give the same \emph{framing} if
$\iota_2^{-1}\circ\iota_1:S^N\rightarrow S^N$ can be deformed inside
the group of diffeomorphisms of $S^N$ to an element of
$O(N+1)$. Identifications that give the same framing define Dehn
twists that are symplectically isotopic. However, it is
unknown whether the Dehn twist can be defined as a symplectic isotopy
class independent of framing. We omit the choice of framing in our
notation as there is often a natural choice in the situations we describe below.
\qed
\end{remark}

Now consider $\CP^n\times\CP^{n+1}$ with the product symplectic form
$\sigma_{\mu}=\sigma_{\CP^n}+\mu\sigma_{\CP^{n+1}}$, where $\mu>1$
and $\sigma_{\CP^n},\sigma_{\CP^{n+1}}$ are the normalized
Fubini-Study forms on $\CP^n$ and $\CP^{n+1}$, respectively. As shown
in \cite[Lemma 1]{CS05}, the graph of the complex conjugate of the Hopf map
$H_{\CP^n}:S^{2n+1}\rightarrow \CP^n$ embeds naturally into this manifold as a Lagrangian sphere $L^{\mu}$. 
\begin{lemma}\label{lem:hopf}
If $\mu>1$, then 
\[
L^{\mu}:=\{[\overline{z}_0:\cdots:\overline{z}_n][z_0:\cdots
:z_n:\sqrt{\mu-1}]\in\CP^n\times\CP^{n+1}\mid \vert
z\vert^2=1\}
\]
is a Lagrangian sphere in $(\CP^n\times\CP^{n+1},\sigma_{\mu})$.
\end{lemma}
\proof For $\lambda>0$, let
$B(\lambda)=\{z\in\CC^{n+1}\mid \vert z\vert^2<\lambda^2\}$.
Then the map $i:B(\sqrt{\mu})\hookrightarrow \CP^{n+1}$ given by
\[
i:(z_0,\ldots,z_n)\mapsto[z_0:\cdots:z_n:\sqrt{\mu-\vert
z_0\vert^2-\cdots-\vert z_n\vert^2}]
\]
is a symplectic embedding
$(B(\sqrt{\mu}),\omega_{\CC^{n+1}})\hookrightarrow
(\CP^{n+1},\mu\sigma_{\CP^{n+1}})$. This follows from the fact that $i$
factors through the Hopf map
$H_{\CP^{n+1}}:S^{2n+3}\rightarrow\CP^{n+1}$
\[
\xymatrix{& S^{2n+3}\ar[d]^{H_{\CP^{n+1}}}\\
B(\sqrt{\mu})\ar[ur]^{\tilde{i}}\ar[r]^i&\CP^{n+1}}
\]
where
\[
\tilde{i}(z_0,\ldots,z_n)=\left(\tfrac{z_0}{\sqrt{\mu}},\ldots,\tfrac{z_n}{\sqrt{\mu}},\sqrt{1-\vert\tfrac{z}{\sqrt{\mu}}\vert^2}\right).
\] 
Hence,
\[
i^*(\mu\sigma_{\CP^{n+1}})=\tilde{i}^*((H_{\CP^{n+1}})^*(\mu\sigma_{\CP^{n+1}}))=\tilde{i}^*(\mu\omega_{\CC^{n+2}}\vert_{S^{2n+3}})=\omega_{\CC^{n+1}}\vert_{B(\sqrt{\mu})}.
\]
Thus, $S^{2n+1}:=\partial B(1)\subset
B(\sqrt{\mu})$ embeds symplectically into
$(\CP^{n+1},\mu\sigma_{\CP^{n+1}})$.  Therefore, the graph of the
conjugate of the Hopf map $H_{\CP^n}:S^{2n+1}\rightarrow \CP^n$
embeds symplectically into $(\CP^n\times\CP^{n+1}, \sigma_{\mu}^n)$
as $L^{\mu}$
. Since
$(\overline{H}_{\CP^n})^*(\sigma_{\CP^n})=-\omega_{\CC^{n+1}}\vert_{S^{2n+1}}$,
\[
\sigma_{\mu}\vert_{L^{\mu}}=\omega_{\CC^{n+1}}\vert_{S^{2n+1}}+(\overline{H}_{\CP^n})^*(\sigma_{\CP^n})=0,
\]
i.e., $L^{\mu}$ is Lagrangian. \qed

\subsection{Lefschetz fibrations}
\begin{definition}\label{def:lef}
A \emph{Lefschetz fibration} is a smooth fibration $\pi:E\rightarrow \CP^1$
such that
\begin{itemize}
\item the set $E_{\crit}$ of critical points of $\pi$ is finite and no two critical
  points lie in the same fiber;
\item $E$ is compact and has a closed $2$-form $\Omega$ that restricts
  to a symplectic form on the non-singular fibers;
\item there is a complex structure $J$ defined in a
  neighborhood of each critical point such that $\Omega$ is $J$-K\"{a}hler;
\item if $j$ is the complex structure on $\CP^1$, then $\pi$ is
  $(J,j)$-holomorphic and at each critical
  point the Hessian of $\pi$ is non-degenerate.
\end{itemize}
\end{definition}
Note that by the Morse Lemma, the fourth condition is equivalent to the
condition that for each $z_{\crit}\in E_{\crit}$ we can find a
$J$-holomorphic coordinate chart $\Psi$ in a
neighborhood of $z_{\crit}$ and a chart $\psi$ on $\CP^1$
centered at $\pi(z_{\crit})$ in $\CP^1$ such that
\[
\psi\circ\pi\circ\Psi^{-1}(z)=\sum_{i=1}^N z_i^2.
\]
Such a pair of charts $(\Psi,\psi)$ will be called a \emph{Morse chart}.
In the following we often suppress the choice of $\psi$ in our
notation.

\begin{figure}[!ht]
\begin{center}
\input{vcycle.pstex_t}
\caption{The vanishing cycle\label{fig:vcycle}}
\end{center}
\end{figure}
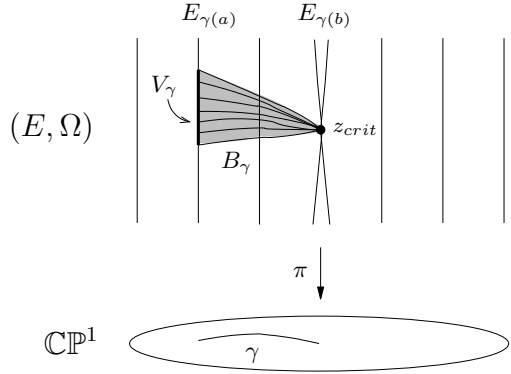

Away from the critical points, $(E,\pi,\Omega)$ is a
symplectic fiber bundle and hence admits a symplectic connection.
Each critical point $z_{\crit}\in E_{\crit}$ gives
rise to a Lagrangian sphere in the non-singular fibers called the
vanishing cycle. Let $\gamma:[a,b]\rightarrow\CP^1$ be an
embedded path which avoids $\pi(E_{\crit})$ except at the endpoint
$\gamma(b)=\pi(z_{\crit})$. Define
\[
B_{\gamma}=\{z_{\crit}\}\cup \bigcup_{a\le s<b}\{z\in E_{\gamma(s)}\mid
\lim_{s'\rightarrow b}\rho_{\gamma\vert_{[s,s']}}(z)=z_{\crit}\}.
\]
By \cite[Lemma 1.13 and 1.14]{Sei03}, $B_{\gamma}$ is an embedded
closed $N$-ball in $E$ with $\Omega\vert_{B_{\gamma}}=0$ whose
boundary
\[
V_{\gamma}=\partial B_{\gamma}=B_{\gamma}\cap E_{\gamma(a)}
\]
is a Lagrangian sphere in
$(E_{\gamma(a)}, \Omega\vert_{E_{\gamma(a)}})$ that comes with a
natural framing (see Remark \ref{framing}). We call $V_{\gamma}$ the
\emph{vanishing cycle} associated to $\gamma$. (See Figure
\ref{fig:vcycle}). If $\gamma'$ is path-homotopic to $\gamma$ in
$\CP^1\setminus\pi(E_{\crit})$, then $V_{\gamma'}$ with its natural
framing is symplectically isotopic to $V_{\gamma}$. Hence, for each
$z_{\crit}\in E_{\crit}$ and $t\in\CP^1\setminus\pi(E_{\crit})$, the
symplectic isotopy class of $V_{\gamma}$ with its natural framing
depends only on the path-homotopy class.

Let $(E,\pi,J,\Omega)$ be a Lefschetz fibration and
$\gamma:[a,b]\rightarrow\CP^1$ an embedded path that avoids
$\pi(E_{\crit})$ except at the endpoint $\gamma(b)=\pi(z_{\crit})\in
\pi(E_{\crit})$. We say that a loop
$\ell:[c,d]\rightarrow\CP^{1}\setminus \pi(E_{\crit})$
\emph{doubles} $\gamma$ if $\ell(c)=\ell(d)=\gamma(a)$, $\ell$ is
positively oriented with respect to the standard orientation of
$\CP^1$ and $\ell$ circles the point $\pi(z_{\crit})$ exactly once and
circles no other critical values of $\pi$. The following result is the
main result of this section and is described in
\cite[Proposition 1.15]{Sei03}. 
\begin{proposition}\label{prop:lef-dehn}
Let $(E,\pi,J,\Omega)$ be a Lefschetz fibration and
$V_{\gamma}\subset E_{t_0}$ the vanishing cycle of a path
$\gamma:[a,b]\rightarrow \CP^1$ with $\gamma(a)=t_0$ and
$\gamma(b)=\pi(z_{\crit})\in E_{\crit}$. Let $\ell$ be a loop that
doubles $\gamma$. The monodromy $\rho_{\ell}:E_{t_0}\rightarrow E_{t_0}$
around $\ell$ is symplectically isotopic to the Dehn twist $\tau_{V_{\gamma}}$
along the vanishing cycle $V_{\gamma}$.
\end{proposition}
\emph{Sketch proof.}
We deform the form $\Omega$ to a closed $2$-form $\Omega^1$ which
is standard with respect to the Morse coordinates at each critical
point. We can choose this deformation such that the monodromy of the
resulting Lefschetz fibration agrees with that of $(E,\pi,J,\Omega)$ up to symplectic isotopy. In a neighborhood of each critical point,
the Lefschetz fibration $(E,\pi,J,\Omega^1)$ agrees with a
neighborhood of $0$ in the standard model
$\pi_{\std}:\CC^N\rightarrow\CC$ described in Section \ref{sec:model}.
Here the vanishing cycles correspond to the $\Sigma_t$'s defined in
(\ref{Sigma_t}). We can therefore use Lemma \ref{lem:stdmodel} to show
the result. Details can be found in \cite[Proposition 2.2.3]{ED}.
\qed

\section{The Lefschetz fibration
  $\pi:\XX\rightarrow\CP^1$}\label{chap:lef}
This section contains the essential steps in the proofs of the main
results of this paper. In the first two sections we construct the
 Lefschetz fibration $(\XX,\pi,J,\Omega^{\mu})$ of Theorem
 \ref{thm:lef}. Its total
space $\XX$ is a complex hypersurface in
$\FF\times\CP^1$, where $\FF$ is a $\CP^{n+1}$-bundle over $\CP^{n+1}$. The
manifold $\XX$ fibers over $\CP^1$ as a subbundle of the
bundle $\FF\times\CP^1\rightarrow\CP^1$. For each $\mu>1$,
the closed $2$-form $\Omega^{\mu}$ on $\XX$ comes from a
toric symplectic structure on $\FF$, and the symplectic isotopy class of the
non-singular fibers of $(\XX,\pi,J,\Omega^{\mu})$ depends
on the cohomology class of this toric structure. For each $\mu>1$, it
is possible to choose a toric structure $\omega_{\FF}^{(1,\mu)}$ on $\FF$ such that the
non-singular fibers of $(\XX,\pi,J,\Omega^{\mu})$  are
isotopic to $(\CP^n\times\CP^{n+1},\sigma_{\mu})$. 
By symplectically embedding a
large part of $(\XX, \Omega^{\mu})$ into the toric manifold
$(\FF,\omega_{\FF}^{(1,\mu)})$, and using a Darboux chart on $\FF$, we
show in Section \ref{sec:vcycle} that the vanishing cycle
$\VV_{\infty}$ in the fiber at infinity is the Lagrangian $L^{\mu}$
defined in (\ref{Lmu}).

\subsection{The toric manifold $(\FF,\omega_{\FF}^{(1,\mu)})$}\label{sec:FF}
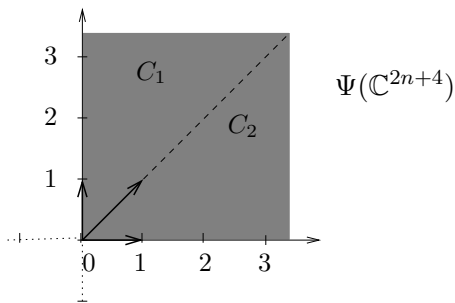
\begin{figure}[!ht]
\begin{center}
\input{FFchambers3.pstex_t}
\caption{The image of $\Psi$ \label{fig:FFchambers3}}
\end{center}
\end{figure}

Let $\CC^{2n+4}$ have coordinates $(\s,q,\x
,x_{n+1})=(s_0,\ldots,s_n,q,x_0,\cdots,x_{n+1})$ and consider the
$(\CC^*)^2$ action on $\CC^{2n+4}$ with weights
\begin{equation}\label{matrix}
\Bigg(\overbrace{\begin{array}{ccc}1&\cdots&1\\0&\cdots&0\end{array}}^{n+1}\begin{array}{c}1\\1\end{array}
\overbrace{\begin{array}{ccc}0&\cdots&0\\1&\cdots&1\end{array}}^{n+2}\Bigg).
\end{equation}
We define $\FF$ to be the quotient
\[
\FF=(\CC^{n+2}\setminus\{0\}\times\CC^{n+2}\setminus\{0\})/(\CC^*)^2.
\]
Since the action of $(\CC^*)^2$ on
$\CC^{n+2}\setminus\{0\}\times\CC^{n+2}\setminus\{0\}$ is free and
proper, this quotient is a smooth complex manifold.\footnote{$\FF$ is
  the GIT quotient $\FF=\CC^{2n+4}\git_{\kappa}(\CC^*)^2$ where
  $\kappa$ is an integral point of the chamber $C_1$ of the effective
  cone depicted in Figure \ref{fig:FFchambers3}. The effective cone is the cone in $\RR^2$
generated by the columns of the matrix (\ref{matrix}). It decomposes naturally into a
union of $2$-dimensional cones, each generated by a pair of column
vectors. The interior of these cones are the chambers $C_1$ and $C_2$. As a complex manifold, the GIT quotient is the
same for different choices of $\kappa$ in the same chamber
\cite[Theorem 3.9]{Tha94}.  }
We
denote the points of $\FF$ by $[[\s,q,\x ,x_{n+1}]]$, where $(\s,q,\x
,x_{n+1})\in\CC^{n+2}\setminus\{0\}\times\CC^{n+2}\setminus\{0\}$
and
\[
[[\s,q,\x ,x_{n+1}]]=[[\alpha \s,\alpha\beta q,\beta\x,\beta x_{n+1}]]\textrm{ for }(\alpha,\beta)\in (\CC^*)^2
\]
The map $\FF\rightarrow\CP^{n+1}$ given by $[[\s,q,\x ,x_{n+1}]]\mapsto[x_0:\cdots:x_{n+1}]$
is clearly well-defined and makes $\FF$ a $\CP^{n+1}$ bundle over
$\CP^{n+1}$.

We now describe a K\"ahler structure on $\FF$ by constructing it as a
symplectic quotient. Let $\CC^{2n+4}$ have the standard symplectic
form $\omega_{\CC^{2n+4}}$ and suppose $T^2$ acts on $\CC^{2n+4}$
with weights (\ref{matrix}). The moment map
$\Psi:\CC^{2n+4}\rightarrow(\mathfrak{t}^2)^*\cong\RR^2$ of the
action is given by
\[
\Psi(\s,q,\x ,x_{n+1})=\left(\sum_{j=0}^n \vert s_j\vert^2+\vert
q\vert^2,\vert q\vert^2+\sum_{j=0}^{n+1}\vert x_j\vert^2\right).
\]
Its image is the cone illustrated in Figure
\ref{fig:FFchambers3} and the regular values are the interiors of the cones $C_1$ and $C_2$. 
By the standard symplectic
quotient construction \cite[Section 5.4]{MS98}, for each $\kappa\in C_1$, the
form $\omega_{\CC^{2n+4}}\mid_{\Psi^{-1}(\kappa)}$ descends to a
symplectic form on the quotient manifold $\Psi^{-1}(\kappa)/T^2$.
The elements of $\Psi^{-1}(\kappa)/T^2$ are denoted $[\s,q,\x
,x_{n+1}]$ where $(\s,q,\x ,x_{n+1})\in\Psi^{-1}(\kappa)$ and
\[
[\s,q,\x ,x_{n+1}]=[e^{i\theta_1}
\s,e^{i(\theta_1+\theta_2)}q,e^{i\theta_2}\x,e^{i\theta_2}x_{n+1}]\textrm{ for }(e^{i\theta_1},e^{i\theta_2})\in T^2.
\]
By \cite[Theorem VII.2.1]{Aud04}, as a complex manifold,
the quotient $\Psi^{-1}(\kappa)/T^2$ equals $\FF$. In fact, the map
$\Psi^{-1}(\kappa)/T^2\rightarrow\FF$
given by
\[
[\s,q,\x,x_{n+1}]\mapsto[[\s,q,\x,x_{n+1}]]
\]
is biholomorphic and the quotient symplectic form is a K\"ahler
structure $\omega_{\FF}^{\kappa}$ on
$\FF$.

\begin{remark}\label{rem:coords}
Points of $\FF$ are denoted $[\s,q,\x ,x_{n+1}]$, when $\FF$ is
thought of as a symplectic quotient $\Psi^{-1}(\kappa)/T^2$ and by
$[[\s,q,\x ,x_{n+1}]]$, when we are interested in the complex
structure on $\FF$.
\qed
\end{remark}

\begin{lemma}\label{lem:FF}
The manifold $\FF$ is a $\CP^{n+1}$ bundle over $\CP^{n+1}$ with
fiber coordinates $[s_0:\cdots:s_n:q]$ and base coordinates
$[x_0:\cdots:x_{n+1}]$. For each $\kappa=(\kappa_1,\kappa_2)\in C_1$, the form $\omega_{\FF}^{\kappa}$ integrates to $\kappa_1\pi$ on a line in the
fiber and to $\kappa_2\pi$ on a line in the
section $s_0=\cdots=s_{n-1}=q=0$. Moreover, the submanifold $\FF\cap(q=0)$ is isomorphic
to the
product $\CP^n\times\CP^{n+1}$ and the restriction $\omega_{\FF}^{\kappa}\mid_{\FF\cap(q=0)}$ is
the product form $\kappa_1\sigma_{\CP^n}+\kappa_2\sigma_{\CP^{n+1}}$.
\end{lemma}
\proof We have already seen the
bundle structure of
$\FF$.

To see the statement about $\FF\cap(q=0)$, note that $\FF\cap(q=0)$
can be thought of as $(\Psi\mid_{(q=0)})^{-1}(\kappa)/T^2$ which
is by construction $\CP^n\times\CP^{n+1}$ with symplectic form
$\kappa_1\sigma_{\CP^n}+\kappa_2\sigma_{\CP^{n+1}}$. More precisely, the
map
\[
(\FF\cap(q=0),\omega_{\FF}^{\kappa}\mid_{(q=0)})\rightarrow
(\CP^n\times\CP^{n+1},\kappa_1\sigma_{\CP^n}+\kappa_2\sigma_{\CP^{n+1}})
\]
given by
\[
[\s,0,\x,x_{n+1}]\mapsto([s_0:\cdots:s_n],[x_0:\cdots:x_{n+1}])
\]
is a biholomorphic symplectomorphism.

The line ($q=s_2=\cdots
=s_n=0$ and $\x=0$) in the fiber ($\x=0$) lies in
$(\FF\cap(q=0),\omega_{\FF}\mid_{(q=0)})$ and corresponds
to
\[
\{([s_0:s_1:0:\cdots:0],[0:\cdots:0:1])\in\CP^n\times\CP^{n+1}\}.
\]
Similarly, the line ($q=s_0=\cdots=s_{n-1}=x_0=\cdots=x_{n-1}$) in the
section ($q=s_0=\cdots=s_{n-1}=0$) corresponds to
\[
\{([0:\cdots:0:1],[0:\cdots:0:x_n:x_{n+1}])\in\CP^n\times\CP^{n+1}\}.
\]
It follows that $\omega_{\FF}^{\kappa}$ integrates to $\kappa_1\pi$
on ($q=s_2=\cdots =s_n=0$ and $\x=0$) and to $\kappa_2\pi$ on
($q=s_0=\cdots=s_{n-1}=x_0=\cdots=x_{n-1}=0$).
 \qed

\begin{figure}[!ht]
\begin{center}
\input{polygrid4.pstex_t}
\caption{Polytope $\Delta^{\kappa}$ representing
   $(\FF,\omega_{\FF}^{\kappa})$ \label{fig:polygrid4} for $n=1$}
\end{center}
\end{figure}
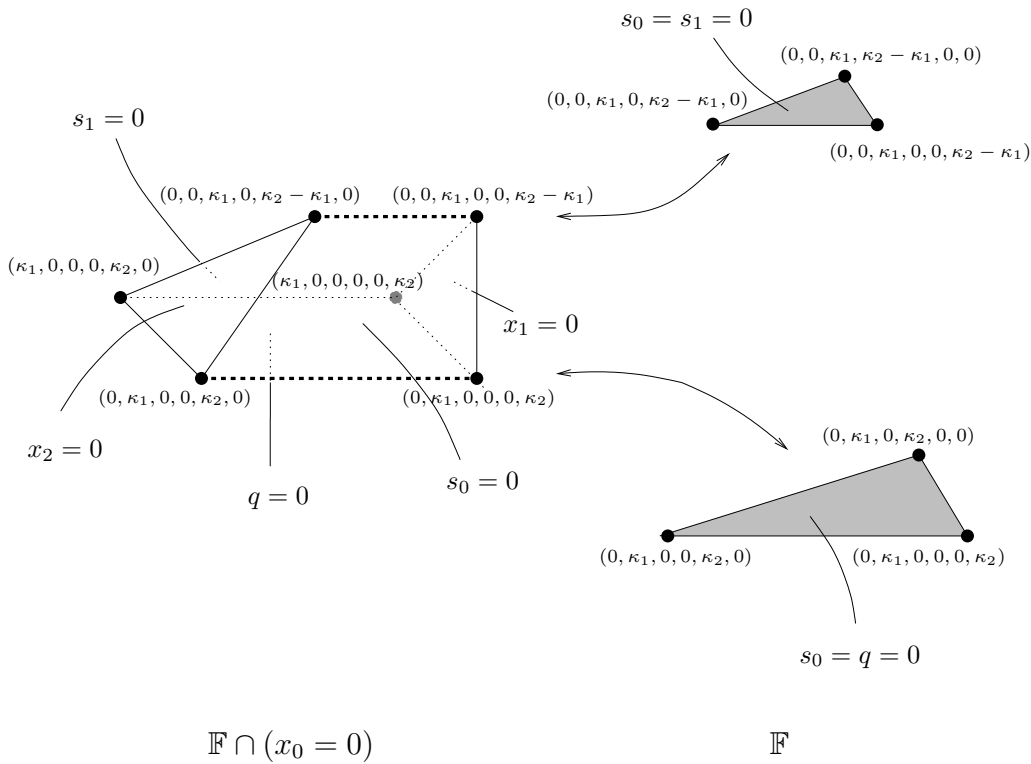

\begin{remark} 
The manifold $(\FF,\omega_{\FF}^{\kappa})$ admits
an action by the quotient torus $T^{2n+4}/T^2\cong T^{2n+2}$. Hence it can be
represented by a Delzant polytope $\Delta^{\kappa}$ of dimension
$2n+2$. If $p:\RR^{2n+4}\rightarrow\RR^2$ is the linear map given by
the matrix (\ref{matrix}), then
\begin{equation}\label{polytope}
\begin{array}{ccl}\Delta^{\kappa}&=&\RR^{2n+4}_+\cap p^{-1}(\kappa)\\
&=&\left\{(\underline{\xi},\nu,\underline{\eta},\eta_{n+1})\in\RR^{2n+4}_+\mid\nu+\sum_{j=0}^n\xi_j=\kappa_1,\:
\nu+\sum_{j=0}^{n+1}\eta_j=\kappa_2\right\}\end{array}
\end{equation}
where $\underline{\xi}=(\xi_0,\ldots,\xi_n)$,
$\underline{\eta}=(\eta_0,\ldots,\eta_n)$ and
$\kappa=(\kappa_1,\kappa_2)$. Let $A:\CC^{2n+4}\rightarrow\RR^{2n+4}$
be the map
\[
A(\s,q,\x,x_{n+1})=(\vert s_0\vert^2,\ldots,\vert s_n\vert^2,\vert q\vert^2,\vert x_0\vert^2,\ldots,\vert x_{n+1}\vert^2).
\]
Then $\Psi^{-1}(\kappa)=A^{-1}(\Delta^{\kappa})$. Suppose we have a
$d$-dimensional face of $\Delta^{\kappa}$. Then it is given by the vanishing of $2n+4-d$ coordinates in
$\RR^{2n+4}$. The vanishing of the corresponding coordinates in
$\CC^{2n+4}$ is a $T^2$-invariant complex submanifold
of $\Psi^{-1}(\kappa)$ and corresponds to a $d$-dimensional
complex submanifold of $\FF=\Psi^{-1}(\kappa)/T^2$. Two faces of
$\Delta^{\kappa}$ intersect precisely when the corresponding
submanifolds intersect in $\FF$.

For $n=1$ the polytope representing $(\FF,\omega_{\FF}^{\kappa})$ is
illustrated in Figure \ref{fig:polygrid4}. The
polytope that represents
$(\FF\cap(x_0=0),\omega_{\FF}^{\kappa}\mid_{x_0=0})$ can be illustrated
in a $3$-dimensional subspace of $\RR^6$ as seen on the left hand
side of Figure \ref{fig:polygrid4}. This is a $\CP^2$-bundle over
$\CP^1$ with base coordinates $[x_1:x_2]$ and
fiber coordinates $[s_0:s_1:q]$. We get the polytope associated
to $(\FF,\omega_{\FF}^{\kappa})$ by replacing each section $\CP^1$ by
a $\CP^2$ as shown in the right hand side of the figure. \qed
\end{remark}

\subsection{The Lefschetz fibration $\pi:\XX\rightarrow\CP^1$}\label{sec:XX}
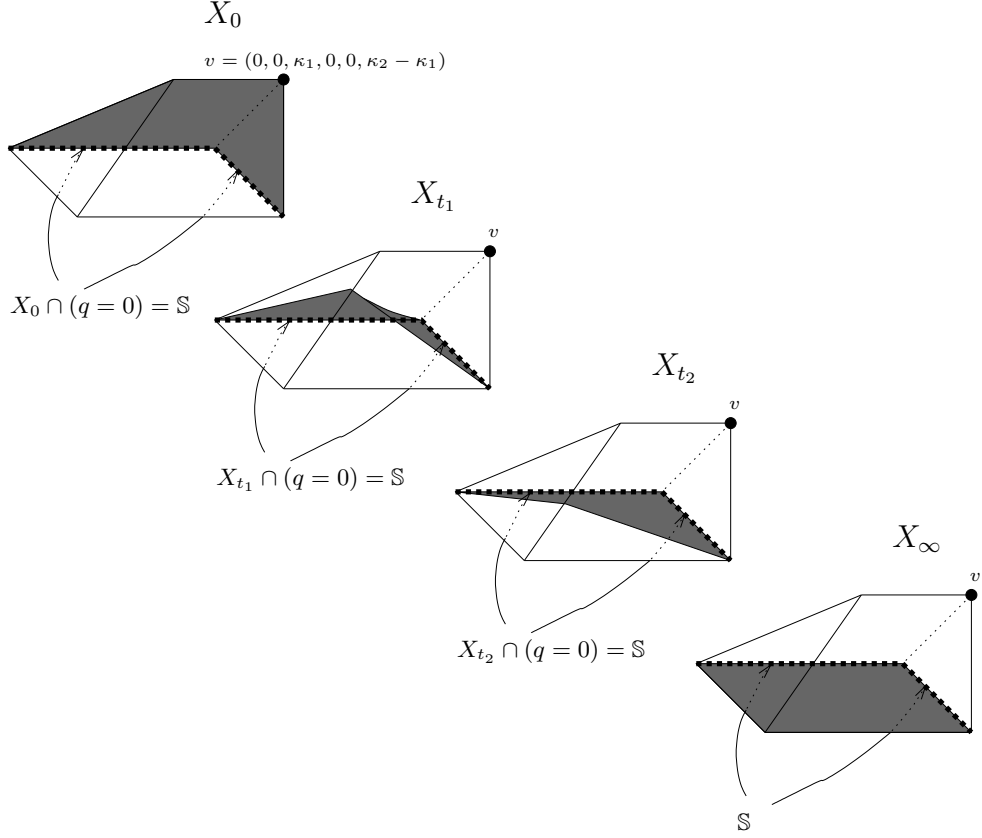
\begin{figure}[!ht]
\begin{center}
\input{fibers2.pstex_t}
\caption{$X_t$ in $\FF\cap(x_0=0)$ \label{fig:fibers2}}
\end{center}
\end{figure}

Let
\[
\XX=\{([[\s,q,\x ,x_{n+1}]],[t_0:t_1])\in\FF\times\CP^1\mid
t_1(s_0x_0+\cdots+s_nx_n)=t_0q\}.
\]
The holomorphic equation $t_1(s_0x_0+\cdots+s_nx_n)=t_0q$ is invariant
under the $(\CC^*)^2$-action on $\CC^{2n+4}$ and $\XX$ is a
complex submanifold of $\FF\times\CP^1$.  Let $t=\tfrac{t_0}{t_1}$
(where we write $\infty$ to mean the point $[1:0]\in\CP^1$) and let
\[
X_t=\{[[\s,q,\x ,x_{n+1}]]\in\FF\mid
t_1(s_0x_0+\cdots+s_nx_n)=t_0q\}
\] 
and $\XX_t=X_t\times\{t\}$. Note that for $t_1=0$ we get
\[
X_{\infty}=\FF\cap (q=0).
\]
Let
$\pi:\XX\rightarrow\CP^1$ be the natural projection map. Figure
\ref{fig:fibers2} shows how the
submanifolds $X_t$ lie inside the polytope representing $(\FF\cap
(x_0=0),\omega_{\FF}^{\kappa}\mid_{x_0=0})$ when $n=1$. 

Consider the submanifold
\begin{equation}\label{def:S}
\SF=\{[[\s,q,\x ,x_{n+1}]]\in\FF\mid
s_0x_0+\cdots+s_nx_n=0\textrm{ and }q=0\}
\end{equation}
in $\FF$. We see that
\[
\SF=\bigcap_{t\in\CP^1}X_t
\]
and moreover, when $\SF$ is removed from the $X_t$'s, they are disjoint in $\FF$ and
\[
\FF=\SF\sqcup\bigsqcup_{t\in\CC}X_t\setminus \SF.
\]
More precisely,
if $\pi_{\FF}:\FF\times\CP^1\rightarrow\FF$ denotes projection to the
first factor, then
\begin{equation}\label{pi_F}
\pi_{\FF}\mid_{\XX\setminus(\SF\times\CP^1)}:\XX\setminus(\SF\times\CP^1)\rightarrow\FF\setminus\SF
\end{equation}
has inverse
$(\pi_{\FF}\mid_{\XX\setminus(\SF\times\CP^1)})^{-1}:\FF\setminus\SF\rightarrow\XX\setminus(\SF\times\CP^1)$ given by
\[
[[\s,q,\x ,x_{n+1}]]\mapsto
\left([[\s,q,\x ,x_{n+1}]],[s_0x_0+\cdots+s_nx_n:q]\right).
\]
In particular, $\pi_{\FF}\mid_{\XX\setminus(\SF\times\CP^1)}$ is biholomorphic.

\begin{lemma}\label{lem:XX}
The singular fibration $\pi:\XX\rightarrow\CP^1$ is holomorphic and has exactly
one critical point at
\[
z_{\crit}=([[\underline{0},1,\underline{0},1]],0)\in
\XX_0.
\]
In a neighborhood of $z_{\crit}$ we obtain holomorphic coordinates
$(z_1,\ldots,z_{2n+2})$ such that $\pi$ is the map
$(z_1,\ldots,z_{2n+2})\mapsto\sum_{j=1}^{2n+2}z_j^2$.

Moreover, when the singular fiber $\XX_0$ is removed,
$\pi:\XX\rightarrow\CP^1$ is trivialized by a
biholomorphic map
$\Phi:\XX\setminus\XX_0\rightarrow(\CP^n\times\CP^{n+1})\times\CC$.
\end{lemma}
\proof Clearly, $\pi:\XX\rightarrow\CC$ is a holomorphic fibration
with exactly one critical point $z_{\crit}$. 

For $q\ne 0$ and $x_{n+1}\ne 0$,
\begin{equation*}
\tilde{s}_i=\frac{s_ix_{n+1}}{q},\quad
\tilde{x}_j=\frac{x_j}{x_{n+1}}\quad \textrm{for }0\le i,j\le n
\end{equation*}
are holomorphic and invariant under the $(\CC^*)^2$-action. The map
\[
\FF\cap(q\ne 0,x_{n+1}\ne 0)\rightarrow \CC^{2n+2}
\] 
given by
\[
[[\s,q,\x ,x_{n+1}]]\mapsto(\tilde{s}_0,\ldots,
\tilde{s}_n,\tilde{x}_0,\ldots,\tilde{x}_n)
\]
defines a holomorphic coordinate chart on $\FF$. Thus, $(\tilde{s}_0,\ldots,
\tilde{s}_n,\tilde{x}_0,\ldots,\tilde{x}_n,t)$ define holomorphic
coordinates on $\FF\times\CC$ in which $\XX$ is given by $\tilde{s}_0\tilde{x}_0+\cdots+\tilde{s}_n\tilde{x}_n=t$. 
Therefore $(\tilde{s}_0,\ldots,
\tilde{s}_n,\tilde{x}_0,\ldots,\tilde{x}_n)$ define holomorphic
coordinates on $\XX$. Let
\begin{equation}
\begin{array}{c}z_1=\tfrac{1}{2}(\tilde{s}_0+\tilde{x}_0),\;\ldots,\;
z_{n+1}=\tfrac{1}{2}(\tilde{s}_n+\tilde{x}_n),\\
z_{n+2}=\tfrac{1}{2i}(\tilde{s}_0-\tilde{x}_0),\;\ldots,\;
z_{2n+2}=\tfrac{1}{2i}(\tilde{s}_n-\tilde{x}_n)\end{array}.
\end{equation}
Then $(z_1,\ldots,z_{2n+2})$ define holomorphic coordinates on $\XX$
in which $\pi$ has the form
$(z_1,\ldots,z_{2n+2})\mapsto\sum_{j=1}^{2n+2}z_j^2$.

Let $\XX_t$ be a non-singular fiber of $\pi:\XX\rightarrow\CP^1$,
i.e., fix $t\ne 0$. Then $X_t$ is a complex submanifold of $\FF$ in
which $q$ is determined by the other coordinates of $\FF$ as
$q=\frac{s_0x_0+\cdots+s_nx_n}{t}$. Since
$\FF=(\CC^{n+2}\setminus\{0\}\times\CC^{n+2}\setminus\{0\})/(\CC^*)^2$,
an element $[[\s,q,\x ,x_{n+1}]]\in X_t$
cannot have $\s=0$ or $\x=x_{n+1}=0$. It follows that the map
$\Phi_t:\XX_t\rightarrow\CP^n\times\CP^{n+1}$ given by
\begin{equation}\label{Phi_t}
\Phi_t:([[\s,q,\x ,x_{n+1}]],t)\mapsto([s_0:\cdots:s_n],[x_0:\cdots:x_{n+1}])
\end{equation}
is well-defined. In fact, it is biholomorphic. Now define
\[
\Phi:\XX\setminus\XX_0\rightarrow
(\CP^n\times\CP^{n+1})\times\CC
\]
by
\begin{equation}\label{Phi}
\Phi([[\s,q,\x ,x_{n+1}]],t)=(\Phi_{\frac{1}{t}}([[\s,q,\x ,x_{n+1}]],\tfrac{1}{t}),\tfrac{1}{t})\quad\textrm{for
}t\in\CP^1\setminus\{0\}.
\end{equation}
\qed
\begin{remark}\label{vertex}
In the notation of Remark
\ref{rem:coords}, 
\[
z_{\crit}=([\underline{0},\kappa_1,\underline{0},\kappa_2-\kappa_1],0).
\]
In the polytope $\Delta^{\kappa}\subset\RR^{2n+4}$
representing $(\FF,\omega_{\FF}^{\kappa})$,
$[\underline{0},\kappa_1,\underline{0},\kappa_2-\kappa_1]\in\FF$ corresponds to the vertex
$v=(\underline{0},\kappa_1,\underline{0},\kappa_2-\kappa_1)$. (See
Figure \ref{fig:fibers2} for the case $n=1$).
\qed
\end{remark}

For each $\mu>1$ we now introduce a $2$-form $\Omega^{\mu}$ on
$\XX$ that makes $(\XX,\pi,J)$ into a
Lefschetz fibration. Here $J$ is the complex structure on
$\XX$ coming from $\FF\times\CP^1$. Recall that if
$\kappa=(\kappa_1,\kappa_2)$ lies in the cone $C_1$, then we have a toric
K\"ahler structure $\omega_{\FF}^{\kappa}$ on $\FF$. If $\mu>1$, then
$(1,\mu)\in C_1$ and we define
\begin{equation}\label{Omega^mu}
\Omega^{\mu}=(\pi_{\FF}\mid_{\XX})^*(\omega_{\FF}^{(1,\mu)}).
\end{equation}
\begin{proposition}\label{prop:isotopy}
For any $\mu>1$, $(\XX,\pi,J,\Omega^{\mu})$ is a
Lefschetz fibration whose non-singular fibers are symplectically
isotopic to $(\CP^n\times\CP^{n+1}, \sigma_{\mu})$. Moreover, the
map
\[
\Phi_{\infty}:(\XX_{\infty},\Omega^{\mu}\mid_{\XX_{\infty}})\rightarrow(\CP^n\times\CP^{n+1},
\sigma_{\mu})\]
described by (\ref{Phi}) is a symplectomorphism.
\end{proposition}
\proof
Let $\mu>1$. First we prove that $(\XX,\pi,J,\Omega^{\mu})$ is a Lefschetz
fibration. Clearly, $\Omega^{\mu}$ is closed. Its restriction
to each non-singular fiber is
\[
\Omega^{\mu}\mid_{\XX_t}=\omega_{\FF}^{(1,\mu)}\mid_{X_t},
\]
which is symplectic since all the $X_t$'s are complex
submanifolds of $\FF$ for $t\ne 0$. We have already shown in Lemma \ref{lem:XX}
that $\pi$ has the form $(z_1,\ldots,z_{2n+2})\mapsto
\sum_{j=1}^{2n+2} z_j^2$ in a neighborhood of the critical point
$z_{\crit}$. By the definition of $\Omega^{\mu}$, the
map
$\pi_{\FF}\mid_{\XX\setminus(\SF\times\CP^1)}$ described in (\ref{pi_F})is a
biholomorphic map that identifies $(\XX\setminus
(\SF\times\CP^1),\Omega^{\mu})$ and $(\FF\setminus
\SF,\omega_{\FF}^{(1,\mu)})$. Hence $\Omega^{\mu}$ is $J$-K\"ahler
in a neighborhood of the critical point. We conclude that
$(\XX,\pi,J,\Omega^{\mu})$ is a Lefschetz fibration.

It remains to show that the non-singular fibers of
$(\XX,\pi,J,\Omega^{\mu})$ are symplectically isotopic to
$(\CP^n\times\CP^{n+1},\sigma_{\mu})$, i.e., that the forms
\[
\alpha_t=(\Phi_t^{-1})^*(\omega_{\FF}^{(1,\mu)}\mid_{X_t})\quad\textrm{for
    }t\ne 0
\]
are symplectically isotopic to $\sigma_{\mu}$ for all $t\ne 0$. Since all the
$\alpha_t$ are symplectically isotopic, it suffices to prove it for $t=\infty$. But
we have already seen in Lemma \ref{lem:FF} that
$\alpha_{\infty}=\sigma_{\mu}$ by construction. \qed

\subsection{The vanishing cycle of $(\XX,\pi,J,\Omega^{\mu})$}\label{sec:vcycle}
In order to compute the vanishing cycle of
$(\XX,\pi,J,\Omega^{\mu})$, we need to understand the
connection coming from the $2$-form $\Omega^{\mu}$ in a neighborhood
of the node $z_{\crit}$. By 
(\ref{Omega^mu}), the map
$\pi_{\FF}\mid_{\XX\setminus(\SF\times\CP^1)}$ described in (\ref{pi_F}) identifies
$(\XX\setminus(\SF\times\CP^1),\Omega^{\mu})$ with
$(\FF\setminus \SF,\omega_{\FF}^{(1,\mu)})$. Hence $\Omega^{\mu}$ is
symplectic on $\XX\setminus(\SF\times\CP^1)$ and we can
use the toric structure of $\FF$ to describe the following Darboux
chart:
\begin{lemma}\label{lem:darboux}
Let
\[
W^{\mu}=\{(\s,\x)\in\CC^{2n+2}\mid \vert \s\vert^2< 1,\; -\vert
\s\vert^2+\vert \x\vert^2< \mu-1\},
\]
where $(\s,\x)=(s_0,\ldots,s_n,x_0,\ldots,x_n)$, $\vert \s\vert^2=\sum_j\vert
  s_j\vert^2$ and $\vert \x\vert^2=\sum_j\vert
  x_j\vert^2$.

The map
\[
\psi_{\XX}^{\mu}:W^{\mu}\rightarrow\XX\cap
(q\ne 0)\cap(x_{n+1}\ne 0)\]
given by
\[
(\s,\x)\mapsto \left(\lbrack \s,\sqrt{1-\vert
  \s\vert^2},\x, \sqrt{\mu-1+\vert \s\vert^2-\vert
  \x\vert^2}\rbrack,\frac{s_0x_0+\cdots+s_nx_n}{\sqrt{1-\vert \s\vert^2}}\right)
\]
defines a Darboux chart on
$(\XX\setminus(\SF\times\CP^1),\Omega^{\mu})$ such that
$\psi_{\XX}^{\mu}(\underline{0},\underline{0})=z_{\crit}$.
\end{lemma}
\proof As noted in the Remark \ref{vertex},
$z_{\crit}=(z,0)\in\FF\times\CP^1$ where 
\[
z=[[\underline{0},1,\underline{0},1]]=[\underline{0},1,\underline{0},\sqrt{\mu-1}].
\]
By the comments
preceding this lemma, it suffices to see that the map  $\psi_{\FF}^{\mu}:W^{\mu}\rightarrow \FF$, where
\[
\psi_{\FF}^{\mu}(\s,\x)=\lbrack \s,\sqrt{1-\vert
  \s\vert^2},\x, \sqrt{\mu-1+\vert \s\vert^2-\vert
  x\vert^2}\rbrack,
\]
is a Darboux chart on $(\FF\cap(q\ne
0)\cap(x_{n+1}\ne 0),\omega_{\FF}^{(1,\mu)})$ centered at $z$.

We have the following commutative diagram
\[
\xymatrix{
{}W^{\mu}\ar[r]^>>>>{\widetilde{\psi}^{\mu}}\ar[dr]_<<<<<<<<<{\psi_{\FF}^{\mu}}&{}\Psi^{-1}((1,\mu))\subset\CC^{2n+4}\ar[d]\\
{}&\FF=\Psi^{-1}((1,\mu))/T^2,}
\]
where $\widetilde{\psi}^{\mu}$ is given by $(\s,\x)\mapsto( \s,\sqrt{1-\vert
  \s\vert^2},\x, \sqrt{\mu-1+\vert \s\vert^2-\vert
  x\vert^2})$. Hence it suffices to see
that  $\widetilde{\psi}^{\mu}$ pulls back $\omega_{\CC^{2n+4}}$ to
$\omega_{\CC^{2n+2}}$. But this is clear since $q,x_{n+1}\in\RR$ on
$\textrm{Im}(\widetilde{\psi}^{\mu})$, and therefore
$dq\wedge d\overline{q}$ and $dx_{n+1}\wedge d\overline{x}_{n+1}$ pull
back to $0$ under $\widetilde{\psi}^{\mu}$.

\qed

\begin{remark}
Each vertex of the polytope $\Delta^{(1,\mu)}$ gives a natural Darboux
chart on $(\FF,\omega_{\FF}^{(1,\mu)})$, and the chart described in Lemma
\ref{lem:darboux} is the one corresponding to $v=(\underline{0},1,\underline{0},\mu-1)$. By (\ref{polytope}), an element $(\underline{\xi},\nu,\underline{\eta},\eta_{n+1})$ of
$\Delta^{(1,\mu)}$ satisfies the equations
\[
\nu=1-\sum_{j=0}^n\xi_j,\quad\eta_{n+1}=\mu-1+\sum_{j=0}^n(\xi_j-\eta_j).
\]
Hence the projection
$(\underline{\xi},\nu,\underline{\eta},\eta_{n+1})\mapsto(\underline{\xi},\underline{\eta})$
maps $\Delta^{(1,\mu)}$ isomorphically to a polytope
$\widetilde{\Delta}^{(1,\mu)}$ in $\RR^{2n+2}$. Here the vertex $v$
is mapped to $0$. Let $F^{q,x_{n+1}}$ denote the union of the facets
$\nu=0$ and $\eta_{n+1}=0$ in $\Delta^{(1,\mu)}$ and let
$\widetilde{F}^{q,x_{n+1}}$ be the corresponding facets in
$\widetilde{\Delta}^{(1,\mu)}$. If
$\widetilde{A}:\CC^{2n+2}\rightarrow\RR^{2n+2}$ is the map
$(\s,\x)\mapsto(\vert s_0\vert^2,\ldots,\vert s_n\vert^2,\vert x_0\vert^2,\ldots,\vert x_n\vert^2)$, then 
$W^{\mu}=\widetilde{A}^{-1}(\widetilde{\Delta}^{(1,\mu)}\setminus\widetilde{F}^{q,x_{n+1}})$
and we have the following commutative diagram:
\[
\xymatrix{
\CC^{2n+2}\supset&W^{\mu}\ar[r]^>>>>>>>>>{\widetilde{\psi}^{\mu}}\ar[d]_{\widetilde{A}}&A^{-1}(\Delta^{(1,\mu)}\setminus
F^{q,x_{n+1}})\ar[d]^{A}&\subset\CC^{2n+4}\\
\RR^{2n+2}\supset&\widetilde{\Delta}^{(1,\mu)}\setminus
\widetilde{F}^{q,x_{n+1}}\ar[r]^{\cong}&\Delta^{(1,\mu)}\setminus F^{q,x_{n+1}}&\subset\RR^{2n+4},
}
\]
where $\widetilde{\psi}^{\mu}$ is the map described in the proof of
Lemma \ref{lem:darboux}.
\qed
\end{remark}

By Lemma \ref{lem:darboux}, we have
a commutative diagram of singular fibrations:
\[
\xymatrix{
(W^{\mu},\omega_{\CC^{2n+2}})\ar[dr]_{\pi_{\loc}}\ar[rr]^(0.4){\psi_{\XX}^{\mu}}_(0.4){\cong}
&&\ar[dl]^{\pi}(\XX\cap(q\ne 0)\cap(x_{n+1}\ne
    0),\Omega^{\mu})\\
 & \CC
}
\]
where 
\[
\pi_{\loc}(\s,\x)=\frac{s_0x_0+\cdots+s_nx_n}{\sqrt{1-\vert
    \s\vert^2}}.
\]
The fiber $W^{\mu}_t=\pi_{\loc}^{-1}(t)$ in $W^{\mu}$ is given by
\[
W^{\mu}_t=\left\{(\s,\x)\in W^{\mu}\mid f_1(\s,\x)=f_2(\s,\x)=0\right\},
\]
where
\begin{eqnarray}
f_1(\s,\x)&=& \sum_{j=0}^n(\Re(s_j)\Re(x_j)-\Im(s_j)\Im(x_j))-\Re(t)\sqrt{1-\vert \s\vert^2},\nonumber\\
f_2(\s,\x)&=&
\sum_{j=0}^n(\Im(s_j)\Re(x_j)+\Re(s_j)\Im(x_j))-\Im(t)\sqrt{1-\vert
\s\vert^2}\nonumber.
\end{eqnarray}
\begin{lemma}\label{lem:connection}
Consider the fibration
$\pi_{\loc}=\pi\circ\psi_{\XX}^{\mu}:W^{\mu}\rightarrow\CC$
with connection coming from the standard symplectic form
$\omega_{\CC^{2n+2}}$ on $W^{\mu}$. Let $t>0$. If
$(\s,\x)\in\pi^{-1}(t)$, then $(\overline{\x},\overline{\s})$ lies in
the horizontal
space $(T_{(\s,\x)}W_t^{\mu})^{\omega_{\RR^{4n+4}}}$.
\end{lemma}
\proof The horizontal space at $(\s,\x)\in W^{\mu}_t$ is the
symplectic complement of $T_{(\s,\x)}W^{\mu}_t$ in
$(\CC^{2n+2},\omega_{\CC^{2n+2}})\cong
(\RR^{4n+4},\omega_{\RR^{4n+4}})$. Note that
\[
T_{(\s,\x)}W^{\mu}_t=(\Span \nabla f_1(\s,\x))^{\perp}\cap(\Span\nabla
f_2(\s,\x))^{\perp}
\]
which has symplectic complement
\begin{eqnarray}
(T_{(\s,\x)}W^{\mu}_t)^{\omega_{\RR^{4n+4}}}&=&\Span\{J_0\nabla
f_1(\s,\x),J_0\nabla f_2(\s,\x)\}.\nonumber
\end{eqnarray}
 Hence the
horizontal space is
spanned by the vectors
\[
\left(\begin{array}{c} \Im(x_0)-\Re(t)\frac{\Im(s_0)}{\sqrt{1-\vert
\s\vert^2}}\\\Re(x_0)+\Re(t)\frac{\Re(s_0)}{\sqrt{1-\vert
\s\vert^2}}\\
\vdots
\\\Im(x_n)-\Re(t)\frac{\Im(s_n)}{\sqrt{1-\vert
\s\vert^2}}\\\Re(x_n)+\Re(t)\frac{\Re(s_n)}{\sqrt{1-\vert
\s\vert^2}}\\\Im(s_0)\\\Re(s_0)\\\vdots\\\Im(s_n)\\\Re(s_n)
\end{array}\right),
\left(\begin{array}{c} -\Re(x_0)-\Im(t)\frac{\Im(s_0)}{\sqrt{1-\vert
\s\vert^2}}\\\Im(x_0)+\Im(t)\frac{\Re(s_0)}{\sqrt{1-\vert
\s\vert^2}}\\\vdots\\-\Re(x_n)-\Im(t)\frac{\Im(s_n)}{\sqrt{1-\vert
\s\vert^2}}\\\Im(x_n)+\Im(t)\frac{\Re(s_n)}{\sqrt{1-\vert
\s\vert^2}}\\-\Re(s_0)\\\Im(s_0)\\\vdots\\-\Re(s_n)\\\Im(s_n)
\end{array}\right).
\]
\qed

Using this, we are able to compute the vanishing cycle of $(\XX,\pi,J,\Omega^{\mu})$.
\begin{proposition}\label{prop:vcycle}
Let $t>0$. The vanishing cycle of
$(\XX,\pi,J,\Omega^{\mu})$ in the fiber $\XX_t$ is isotopic
to
\[
\VV^{\mu}_t=\left\{([\overline{\x},\sqrt{1-g(t)},\x,\sqrt{\mu-1}],t)\in\XX\mid
\vert \x\vert^2=g(t)\right\}
\]
where $\overline{\x}=(\overline{x}_0,\ldots,\overline{x}_n)$ and
\begin{eqnarray}
g(t)=\frac{-t^2+\sqrt{t^4+4t^2}}{2}.
\end{eqnarray}
The vanishing cycle in $\XX_{\infty}$ is isotopic to
\[
\VV^{\mu}_{\infty}=\left\{([\overline{\x},0,\x,\sqrt{\mu-1}],\infty)\in\XX\mid
\vert \x\vert^2=1\right\}.
\]
\end{proposition}
\proof
For each $t>0$, define
\[
L^{\mu}_t=\left\{(\overline{\x},\x)\in W^{\mu}\mid \vert
\x\vert^2=g(t)\right\}.
\]
The condition on $\vert \x\vert^2$ ensures that $L^{\mu}_t\subset
W^{\mu}_t$. We claim that parallel transport along the path
$\gamma(r)=r$ for $0<r\le t$ takes all of $L_t^{\mu}$ to the node
$(\underline{0},\underline{0})\in W^{\mu}_0$. Since $L_t^{\mu}$ is
diffeomorphic to $S^{2n+1}$ and 
\[
\psi_{\overline{\XX}}^{\mu}(L^{\mu}_t)=\VV^{\mu}_t,
\]
this shows that $\VV^{\mu}_t$ is the vanishing cycle in $\XX_t$.

Fix $t>0$ and let $\gamma_t(r)=r$ for $0<r\le t$. Given
$(\overline{\x},\x)\in L_t^{\mu}$, define the path
\[
\widetilde{\gamma_t}(r)=h_t(r)(\overline{\x},\x),\quad 0<r\le t
\]
where
\[
h_t(r)=\sqrt{\frac{g(r)}{g(t)}}.
\]
This is a lift of $\gamma_t$ with endpoints
$(\overline{\x},\x)$ and $(\underline{0},\underline{0})$ and tangents
\[
\frac{d}{dr}\widetilde{\gamma_t}(r)=h_t'(r)(\overline{\x},\x).
\]
For all $r$, $\widetilde{\gamma_t}(r)\in L_{r}^{\mu}$, and
by Lemma \ref{lem:connection}, the horizonal space at
$\widetilde{\gamma_t}(r)$ contains $\frac{d}{dr}\widetilde{\gamma_t}(r)$. Hence parallel transport along
$\gamma_t$ takes $L_t^{\mu}$ to the node. Since
$\lim_{t\rightarrow\infty}g(t)=1$,
\[
\VV_{\infty}^{\mu}=\VV_{\gamma_{\infty}}^{\mu}=\lim_{t\rightarrow\infty}\VV_{\gamma_t}^{\mu}=\left\{([\overline{\x},0,\x,\sqrt{\mu-1}],\infty)\in\XX\mid
\vert \x\vert^2=1\right\}.
\]
\qed
\section{Proofs of the main results}
We are now in a position to prove the main results described in
Section \ref{chap:intro}.

\medskip

\noindent \emph{Proof of Theorem \ref{thm:lef}.} Given $\mu>1$,
Proposition \ref{prop:isotopy} shows that
$(\XX,\pi,J,\Omega^{\mu})$ is a Lefschetz fibration with
one critical point
\[
z_{\crit}=([\underline{0},1,\underline{0},\sqrt{\mu-1}],0)\in\XX_{0}.
\]
The biholomorphic trivialization
$\Phi:\XX\setminus\XX_0:\rightarrow(\CP^n\times\CP^{n+1})\times\CC$
defined by (\ref{Phi}) induces a symplectomorphism
\[
\Phi_{\infty}:(\XX_{\infty},\Omega^{\mu}\mid_{\XX_{\infty}})\cong(\CP^n\times\CP^{n+1},\sigma_{\mu})
\]
and by Proposition \ref{prop:vcycle}, the vanishing cycle in
$\XX_{\infty}$ is isotopic to
\[
\VV^{\mu}_{\infty}=\left\{([\overline{\x},0,\x,\sqrt{\mu-1}],\infty)\in\XX\mid
\vert \x\vert^2=1\right\}.
\]
Under the symplectomorphism
$\Phi_{\infty}$ this corresponds to
\[
\Phi_{\infty}(\VV_{\infty}^{\mu})=\{([\overline{x}_0:\cdots:\overline{x}_n],[x_0:\cdots:x_n:\sqrt{\mu-1}])\in\CP^n\times\CP^{n+1}\mid
\vert \x\vert^2=1\},
\]
which is exactly the Lagrangian $L^{\mu}$ in $\CP^n\times\CP^{n+1}$.
\qed

\medskip

\noindent \emph{Proof of Theorem \ref{thm:isotopy}}.
Given $\mu>1$, we consider the Lefschetz fibration
$(\XX,\pi,J,\Omega^{\mu})$ from Theorem \ref{thm:lef}.
By Proposition \ref{prop:lef-dehn}, the monodromy
$\rho_{\ell}:\XX_{\infty}\rightarrow\XX_{\infty}$ along a positively
oriented loop $\ell$ in $\CP^1$
based at $\infty\in\CP^1$ and
circling $0$ exactly once is symplectically isotopic to the Dehn twist
$\tau_{\VV_{\infty}^{\mu}}$ along
the vanishing cycle $\VV_{\infty}^{\mu}$.
Since the fibration has only
one singular fiber, $\rho_{\ell}$
is symplectically isotopic to the identity through an isotopy that
arises from deforming the loop $\ell$ to the constant loop based at
infinity in $\CP^1$. Hence we have an isotopy $(\varphi_{\lambda})_{0\le
  \lambda\le 1}$ of $\XX_{\infty}$ such that $\varphi_1$ equals $\tau_{\VV_{\infty}}$ and $\varphi_0$ is the identity.
Recall that each fiber of $\XX$
contains the hypersurface $\SF\times\{t\}$, with $\SF$ defined in
(\ref{def:S}). Parallel transport in
$(\XX,\pi,J,\Omega^{\mu})$ along a path in $\CP^1$ from $t'$ to
$t$ takes any point
$(z,t)\in\SF\times\{t\}\subset\XX_{t}$ to
$(z,t')\in\SF\times\{t'\}\subset\XX_{t'}$.
Hence $\varphi_{\lambda}$ restricts to the identity
on $\SF\times\{\infty\}$ for all $\lambda$. Observe that
$S=\Phi_{\infty}(\SF\times\{\infty\})$ and that by property \ref{lefdelta3} of
Theorem \ref{thm:lef}
\[
\tau_{L^{\mu}}=\Phi_{\infty}\circ\tau_{\VV_{\infty}^{\mu}}\circ\Phi_{\infty}^{-1}.
\]
Hence $\Phi_{\infty}\circ\varphi_{\lambda}\circ\Phi_{\infty}^{-1}$, $0\le
\lambda\le 1$, is the desired isotopy of $\tau_{L^{\mu}}$.
\qed

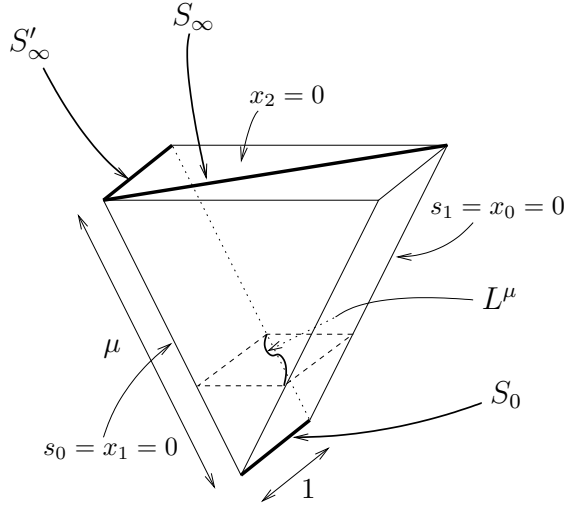
\begin{figure}[!ht]
\begin{center}
\input{P1xP2.pstex_t}
\caption{Polytope representing $(\CP^1\times\CP^2,\sigma_{\mu})$ \label{fig:P1xP2}}
\end{center}
\end{figure}

Under the projection to the first factor, the complex hypersurface $S$ in
$\CP^n\times\CP^{n+1}$ defined in Theorem \ref{thm:isotopy} is a $\CP^n$-bundle over $\CP^n$. The fibers are linearly embedded copies of
$\CP^n$ in $\CP^{n+1}$. As before, let
\[
\begin{array}{ccl}
S_0&=&\CP^n\times\{[0:\cdots:0:1]\}\\
S_{\infty}&=&\{([s_0:\cdots:s_n],[x_0:x_1:0:\cdots:0])\mid
s_0x_0+s_1x_1=0\};\end{array}
\]
these are sections of $S$. If $n=1$, $S$ is a Hirzebruch surface, $S_0$ is the small section (with negative self-intersection), and
$S_{\infty}$ is the large section (with positive
self-intersection). As stated in Corollary \ref{cor:blow-up}, if we
symplectically blow up $\CP^n\times\CP^{n+1}$ along $S_0$
and $S_{\infty}$, then as long as the
blow-ups are small,  the isotopy of the Dehn twist
$\tau_{L^{\mu}}$ remains. On the other hand, Proposition
\ref{prop:cannotfix} (which is a consequence of the proof of
\cite[Proposition 2]{CS05}) shows that an isotopy cannot fix the
submanifolds $S_0$ and $S_{\infty}'$.
Figure \ref{fig:P1xP2} illustrates how the submanifolds $S_0$,
$S_{\infty}$ and $S_{\infty}'$ lie inside the polytope
representing $(\CP^1\times\CP^2,\sigma_{\mu})$.
 
Before proving Corollary
\ref{cor:blow-up} and Proposition \ref{prop:cannotfix}, we review the
symplectic blow-up construction as described in \cite[Chapter
7]{MS98}. First we consider a blow-up at the origin in
$\CC^N$.
Let $\mathcal{L}$ be the tautological line bundle over $\CP^{N-1}$,
i.e., the total space of $\mathcal{L}$ is given by
\[
\mathcal{L}=\{((z_1,\ldots,z_N),[w_1,\ldots,w_N])\in \CC^N\times\CP^{N-1}\mid
w_jz_k=w_kz_j\textrm{ for all }j,k\}.
\]
Let $\pr_{\CC^N}:\mathcal{L}\rightarrow
\CC^N$ and $\pr_{\CP^{N-1}}:\mathcal{L}\rightarrow\CP^{N-1}$ be the
natural projections. For each
$\lambda>0$ the $2$-form
\[
\omega_{\mathcal{L}}^{\lambda}=\pr_{\CC^N}^*(\omega_{\CC^N})+\lambda^2\pr_{\CP^{N-1}}^*(\sigma_{\CP^{N-1}})
\]
is symplectic.
For all $\delta\ge 0$, let $\mathcal{L}(\delta)=\pr_{\CC^N}^{-1}(B(\delta))$,
where $B(\delta)$ is the ball of radius $\delta$ in $\CC^N$.
By \cite[Lemma 7.11]{MS98}, for all $\lambda,\delta>0$
$(\mathcal{L}(\delta)-\mathcal{L}(0),\omega_{\mathcal{L}}^{\lambda})$ is symplectomorphic to the
spherical shell
$(B(\sqrt{\lambda^2+\delta^2})-B(\lambda),\omega_{\CC^N})$. The
blow-up $\widetilde{\CC}^N$ of $0$ in $\CC^N$ of size $\lambda$ is obtained by replacing
the interior of the ball $B(\sqrt{\lambda^2+\delta^2})$ in $\CC^N$
with $\mathcal{L}(\delta)$. Thus, we have replaced the ball $B(\lambda)$ in $\CC^N$ with
the zero section $\mathcal{L}(0)\cong\CP^{N-1}$, the exceptional divisor in
$\widetilde{\CC}^N$. The manifold $\widetilde{\CC}^N$ has a natural
symplectic structure that is independent of the choice of
$\delta>0$. 


Let $Q$ be a compact symplectic submanifold of codimension $2N$ in a
symplectic manifold $(M,\omega)$. Then $\omega$ can be written in a standard form in a
neighborhood of $Q$. Details can be found in \cite[p.
249-251]{MS98}.  The normal bundle $\nu_{Q\mid M}$ of $Q$ in $M$ has
the structure of a $2N$-dimensional symplectic vector bundle with a compatible complex
  structure $J$. It is
  associated to a principal $U(N)$-bundle $P$, i.e.,  $\nu_{Q\mid
  M}=P\times_{U(N)}\CC^N$. This bundle has a symplectic form
  \begin{equation}\label{Qstandardform}
\omega'=\alpha+\pi_Q^*(\omega\mid_{Q}),
\end{equation}
where $\pi_Q:\nu_{Q\mid
    M}\rightarrow Q$ is the natural projection and $\alpha$
    is a closed $2$-form on $\nu_{Q\mid M}$ that restricts
    to the standard symplectic form on the fibers $\CC^N$. A
  small neighborhood of $Q$ in $(M,\omega)$
  is symplectomorphic to the disk bundle $\nu_{Q\mid
    M}^{<\varepsilon}=P\times_{U(N)}B(\varepsilon)$ with symplectic
  form $\omega'$ for some $\varepsilon>0$.
 Using the
standard form of the symplectic structure near $Q$, we can define
small blow-ups along $Q$ in $M$. If $\lambda<\varepsilon$, then the
size $\lambda$ blow-up of $Q$ in $M$ is simply the result of blowing
up the origin in the fiber $B(\varepsilon)$, i.e. if
$\widetilde{B}(\varepsilon)$ is the size $\lambda$ blow-up of the
origin in $B(\varepsilon)$, then the size $\lambda$ blow-up of $Q$
in $M$ is obtained by removing a neighborhood of $Q$ and glueing in
the manifold $P\times_{U(N)}\widetilde{B}(\varepsilon)$ in the
obvious way.

\medskip

\noindent\emph{Proof of Corollary \ref{cor:blow-up}.}
Under the trivialization $\Phi$ of
$(\XX,\pi,J,\Omega^{\mu})$ defined in (\ref{Phi}), the section $S_0$ in $\CP^n\times\CP^{n+1}$ corresponds to
\[
\SF_0=\{[\s,q,\x ,x_{n+1}]\in\FF\mid q=0,\;\x=\underline{0}\}.
\]
Let $Q=\SF_0\times\CP^1$. 

For all $t\in\CP^1$, the fiber $\XX_t$ contains $\SF_0\times\{t\}$ as a
symplectic submanifold and
the diagram
\[
\xymatrix{\XX_t\ar[r]&\XX\ar[r]&\CP^1\\
\SF_0\times\{t\}\ar[r]\ar@{^{(}->}[u]&Q\ar[r]\ar@{^{(}->}[u]&\CP^1\ar@{=}[u]}
\]
commutes.
It follows that the normal bundle of $Q$ in $\XX$ can be
obtained by glueing the normal bundles of $\SF_0\cong\SF_0\times\{t\}$
in $X_t\cong\XX_t$ for all $t\in\CP^1$. Now let
$\widetilde{\XX}$ denote a small symplectic blow-up along
$Q$ in $\XX$ in the neighborhood of $Q$ on which
$(\Omega^{\mu})'$ has the standard form (\ref{Qstandardform}). Then each fiber in
$\widetilde{\XX}$ is of the form
$\widetilde{X}_t\times\{t\}$, where
$\widetilde{X}_t$ is the blow-up of $\SF_0$ in $X_t$. Hence the
isotopy of $\tau_{L^{\mu}}$ lifts to the blow-up of $S_0$ in $\CP^n\times\CP^{n+1}$.

To see that we can blow up along $S_{\infty}$ in $\CP^n\times\CP^{n+1}$, we
repeat this construction with $\SF_0$ replaced by
\[
\SF_{\infty}=\{[\s,0,\x,0]\in\FF\mid x_j=0 \textrm{ for }j\in\{2,\ldots,n\},\;s_0x_0+s_1x_1=0\}
\]
and $Q$ by the complex submanifold $\SF_{\infty}\times\CP^1$ in $\XX$.
\qed

\medskip
\noindent\emph{Proof of Proposition \ref{prop:cannotfix}.} Assume we
have an isotopy $(\varphi_{\lambda})_{0\le\lambda\le 1}$ such that
$\varphi_1=\tau_{L^{\mu}}$ and $\varphi_{\lambda}$ restricts to the
identity on $S_0\sqcup S_{\infty}'$ for all $\lambda$. We can then
deform $(\varphi_{\lambda})_{0\le\lambda\le
  1}$ to a homotopy $(\widetilde{\varphi}_{\lambda})_{0\le\lambda\le
  1}$ from
$\tau_{L^{\mu}}$ to the identity such that each $\widetilde{\varphi}_{\lambda}$ fixes the neighborhoods
\[
\begin{array}{ccccl}
\mathcal{N}(S_0)&=&\CP^n\times B(\varepsilon,p_0)&\supset& S_0,\\
\mathcal{N}(S_{\infty}')&=&\CP^n\times B(\varepsilon,p_{\infty})&\supset& S_{\infty}',
\end{array}
\]
where $B(\varepsilon,p)$ denotes an embedded size $\varepsilon$
(real) $(2n+2)$-ball in $\CP^{n+1}$ centered at $p=p_0,p_{\infty}$. We show that this leads to a contradiction as
follows. Following \cite[p. 534-535]{Gom95}, we create a new smooth manifold $M$ by performing a smooth
surgery on $\CP^n\times\CP^{n+1}$ that takes place in the
neighborhoods $\mathcal{N}(S_0)$ and $\mathcal{N}(S_{\infty})$. In the manifold $M$, the Dehn twist
corresponds to a diffeomorphism that acts non-trivially on homology
but is also homotopic to the identity. This is a contradiction.

Choose a diffeomorphism
$\phi:B(\tfrac{\varepsilon}{2})\setminus\{0\}\rightarrow
B(\tfrac{\varepsilon}{2})\setminus\{0\}$ of the punctured
$(2n+2)$-ball that turns it smoothly inside out, e.g.,
\[
\phi: z\mapsto z\sqrt{\tfrac{\varepsilon^2}{4\vert z\vert^2}-1}.
\]
Let
$\widehat{\CP}^{n+1}$ be the manifold obtained by removing the points
$p_0$ and $p_{\infty}$ from $\CP^{n+1}$ and identifying the open sets
$B(\tfrac{\varepsilon}{2},p_0)\setminus\{p_0\}$ and $B(\tfrac{\varepsilon}{2},p_{\infty})\setminus\{p_{\infty}\}$ via the
map $\phi$. Let
\[
M=\CP^n\times\widehat{\CP}^{n+1}.
\]
Since the surgery took place on neighborhoods on which the Dehn
twist $\tau_L$ and the homotopy $(\widetilde{\varphi}_{\lambda})_{0\le\lambda\le 1}$ restrict to the
identity, $\tau_L$ induces a diffeomorphism $\widehat{\tau}_L$ of $M$
homotopic to the
identity. The manifold $M$ contains a $(2n+1)$-cycle $C$ with
non-trivial intersection with $L$. Indeed, let
\[
\Gamma=\{[x_0:0:\cdots:0:\sqrt{\mu-x_0}]\in\CP^{n+1}\mid 0\le x_0\le\mu\}
\]
be the image of a path from $p_0$ to $p_{\infty}$ in $\CP^{n+1}$. Then
$\CP^n\times\Gamma$ gives rise to a $(2n+1)$-cycle $C$ in $M$ that
intersects the Lagrangian $(2n+1)$-sphere $L^{\mu}$ exactly once,
namely at the point
\[
([1:0\cdots:0],[1:0:\cdots:0:\sqrt{\mu-1}])\in M,
\]
and the intersection is transverse.
In particular, $L$ represents a non-trivial element in the homology
of $M$. By \cite[(1.5)]{ST01}, the Dehn twist $\widehat{\tau}_L$
acts on the homology class of $C$ by adding the
homology class of $L$ to it. See also Figure \ref{fig:dehntwist}.
\qed

\bibliography{sympisotopy}

\end{document}

%% file: dehntwistuv.pstex_t
\begin{picture}(0,0)%
\includegraphics{dehntwistuv.pstex}%
\end{picture}%
\setlength{\unitlength}{2487sp}%
\begingroup\makeatletter\ifx\SetFigFont\undefined%
\gdef\SetFigFont#1#2#3#4#5{%
  \reset@font\fontsize{#1}{#2pt}%
  \fontfamily{#3}\fontseries{#4}\fontshape{#5}%
  \selectfont}%
\fi\endgroup%
\begin{picture}(5471,4515)(2389,-4723)
\put(6226,-376){\makebox(0,0)[lb]{\smash{{\SetFigFont{9}{10.8}{\familydefault}{\mddefault}{\updefault}{\color[rgb]{0,0,0}$\tau(F_u)$}%
}}}}
\put(7291,-2191){\makebox(0,0)[lb]{\smash{{\SetFigFont{9}{10.8}{\familydefault}{\mddefault}{\updefault}{\color[rgb]{0,0,0}$L_0$}%
}}}}
\put(2686,-4651){\makebox(0,0)[lb]{\smash{{\SetFigFont{12}{14.4}{\familydefault}{\mddefault}{\updefault}{\color[rgb]{0,0,0}$T^*S^N$}%
}}}}
\put(2611,-376){\makebox(0,0)[lb]{\smash{{\SetFigFont{9}{10.8}{\familydefault}{\mddefault}{\updefault}{\color[rgb]{0,0,0}$F_u\cong\RR^N$}%
}}}}
\put(3676,-2191){\makebox(0,0)[lb]{\smash{{\SetFigFont{9}{10.8}{\familydefault}{\mddefault}{\updefault}{\color[rgb]{0,0,0}$L_0\cong S^N$}%
}}}}
\put(5086,-1741){\makebox(0,0)[lb]{\smash{{\SetFigFont{9}{10.8}{\familydefault}{\mddefault}{\updefault}{\color[rgb]{0,0,0}$(-u,0)$}%
}}}}
\put(4621,-3511){\makebox(0,0)[lb]{\smash{{\SetFigFont{11}{13.2}{\familydefault}{\mddefault}{\updefault}{\color[rgb]{0,0,0}$\tau$}%
}}}}
\put(3025,-2521){\makebox(0,0)[lb]{\smash{{\SetFigFont{9}{10.8}{\familydefault}{\mddefault}{\updefault}{\color[rgb]{0,0,0}$(u,0)$}%
}}}}
\put(6301,-4651){\makebox(0,0)[lb]{\smash{{\SetFigFont{12}{14.4}{\familydefault}{\mddefault}{\updefault}{\color[rgb]{0,0,0}$T^*S^N$}%
}}}}
\end{picture}%

%% file: vcycle.pstex_t
\begin{picture}(0,0)%
\includegraphics{vcycle.pstex}%
\end{picture}%
\setlength{\unitlength}{2526sp}%
\begingroup\makeatletter\ifx\SetFigFont\undefined%
\gdef\SetFigFont#1#2#3#4#5{%
  \reset@font\fontsize{#1}{#2pt}%
  \fontfamily{#3}\fontseries{#4}\fontshape{#5}%
  \selectfont}%
\fi\endgroup%
\begin{picture}(4913,3625)(541,-4838)
\put(3331,-3901){\makebox(0,0)[lb]{\smash{{\SetFigFont{9}{10.8}{\familydefault}{\mddefault}{\updefault}{\color[rgb]{0,0,0}$\pi$}%
}}}}
\put(3721,-2476){\makebox(0,0)[lb]{\smash{{\SetFigFont{8}{9.6}{\familydefault}{\mddefault}{\updefault}{\color[rgb]{0,0,0}$z_{crit}$}%
}}}}
\put(2236,-1381){\makebox(0,0)[lb]{\smash{{\SetFigFont{9}{10.8}{\familydefault}{\mddefault}{\updefault}{\color[rgb]{0,0,0}$E_{\gamma(a)}$}%
}}}}
\put(3376,-1381){\makebox(0,0)[lb]{\smash{{\SetFigFont{9}{10.8}{\familydefault}{\mddefault}{\updefault}{\color[rgb]{0,0,0}$E_{\gamma(b)}$}%
}}}}
\put(2866,-4696){\makebox(0,0)[lb]{\smash{{\SetFigFont{9}{10.8}{\familydefault}{\mddefault}{\updefault}{\color[rgb]{0,0,0}$\gamma$}%
}}}}
\put(541,-2521){\makebox(0,0)[lb]{\smash{{\SetFigFont{12}{14.4}{\familydefault}{\mddefault}{\updefault}{\color[rgb]{0,0,0}$(E,\Omega)$}%
}}}}
\put(886,-4666){\makebox(0,0)[lb]{\smash{{\SetFigFont{12}{14.4}{\familydefault}{\mddefault}{\updefault}{\color[rgb]{0,0,0}$\mathbb{CP}^1$}%
}}}}
\put(2626,-2836){\makebox(0,0)[lb]{\smash{{\SetFigFont{9}{10.8}{\familydefault}{\mddefault}{\updefault}{\color[rgb]{0,0,0}$B_{\gamma}$}%
}}}}
\put(1936,-2086){\makebox(0,0)[lb]{\smash{{\SetFigFont{9}{10.8}{\familydefault}{\mddefault}{\updefault}{\color[rgb]{0,0,0}$V_{\gamma}$}%
}}}}
\end{picture}%

%% file: FFchambers3.pstex_t
\begin{picture}(0,0)%
\includegraphics{FFchambers3.pstex}%
\end{picture}%
\setlength{\unitlength}{2526sp}%
\begingroup\makeatletter\ifx\SetFigFont\undefined%
\gdef\SetFigFont#1#2#3#4#5{%
  \reset@font\fontsize{#1}{#2pt}%
  \fontfamily{#3}\fontseries{#4}\fontshape{#5}%
  \selectfont}%
\fi\endgroup%
\begin{picture}(3237,2979)(2914,-5323)
\put(4201,-3061){\makebox(0,0)[lb]{\smash{{\SetFigFont{10}{12.0}{\familydefault}{\mddefault}{\updefault}{\color[rgb]{0,0,0}$C_1$}%
}}}}
\put(5101,-3586){\makebox(0,0)[lb]{\smash{{\SetFigFont{10}{12.0}{\familydefault}{\mddefault}{\updefault}{\color[rgb]{0,0,0}$C_2$}%
}}}}
\put(6151,-3211){\makebox(0,0)[lb]{\smash{{\SetFigFont{11}{13.2}{\familydefault}{\mddefault}{\updefault}{\color[rgb]{0,0,0}$\Psi(\mathbb{C}^{2n+4})$}%
}}}}
\put(3676,-4936){\makebox(0,0)[lb]{\smash{{\SetFigFont{10}{12.0}{\familydefault}{\mddefault}{\updefault}{\color[rgb]{0,0,0}$0$}%
}}}}
\put(3301,-3511){\makebox(0,0)[lb]{\smash{{\SetFigFont{10}{12.0}{\familydefault}{\mddefault}{\updefault}{\color[rgb]{0,0,0}$2$}%
}}}}
\put(4180,-4936){\makebox(0,0)[lb]{\smash{{\SetFigFont{10}{12.0}{\familydefault}{\mddefault}{\updefault}{\color[rgb]{0,0,0}$1$}%
}}}}
\put(5401,-4936){\makebox(0,0)[lb]{\smash{{\SetFigFont{10}{12.0}{\familydefault}{\mddefault}{\updefault}{\color[rgb]{0,0,0}$3$}%
}}}}
\put(4816,-4936){\makebox(0,0)[lb]{\smash{{\SetFigFont{10}{12.0}{\familydefault}{\mddefault}{\updefault}{\color[rgb]{0,0,0}$2$}%
}}}}
\put(3301,-2911){\makebox(0,0)[lb]{\smash{{\SetFigFont{10}{12.0}{\familydefault}{\mddefault}{\updefault}{\color[rgb]{0,0,0}$3$}%
}}}}
\put(3301,-4111){\makebox(0,0)[lb]{\smash{{\SetFigFont{10}{12.0}{\familydefault}{\mddefault}{\updefault}{\color[rgb]{0,0,0}$1$}%
}}}}
\end{picture}%

%% file: polygrid4.pstex_t
\begin{picture}(0,0)%
\includegraphics{polygrid4.pstex}%
\end{picture}%
\setlength{\unitlength}{2684sp}%
\begingroup\makeatletter\ifx\SetFigFont\undefined%
\gdef\SetFigFont#1#2#3#4#5{%
  \reset@font\fontsize{#1}{#2pt}%
  \fontfamily{#3}\fontseries{#4}\fontshape{#5}%
  \selectfont}%
\fi\endgroup%
\begin{picture}(8930,6971)(1,-6841)
\put(4057,-4327){\makebox(0,0)[lb]{\smash{{\SetFigFont{10}{12.0}{\familydefault}{\mddefault}{\updefault}{\color[rgb]{0,0,0}$s_0=0$}%
}}}}
\put(4582,-2869){\makebox(0,0)[lb]{\smash{{\SetFigFont{10}{12.0}{\familydefault}{\mddefault}{\updefault}{\color[rgb]{0,0,0}$x_1=0$}%
}}}}
\put(188,-4026){\makebox(0,0)[lb]{\smash{{\SetFigFont{10}{12.0}{\familydefault}{\mddefault}{\updefault}{\color[rgb]{0,0,0}$x_2=0$}%
}}}}
\put(7140,-355){\makebox(0,0)[lb]{\smash{{\SetFigFont{7}{8.4}{\familydefault}{\mddefault}{\updefault}{\color[rgb]{0,0,0}$(0,0,\kappa_1,\kappa_2-\kappa_1,0,0)$}%
}}}}
\put(7593,-1252){\makebox(0,0)[lb]{\smash{{\SetFigFont{7}{8.4}{\familydefault}{\mddefault}{\updefault}{\color[rgb]{0,0,0}$(0,0,\kappa_1,0,0,\kappa_2-\kappa_1)$}%
}}}}
\put(602,-968){\makebox(0,0)[lb]{\smash{{\SetFigFont{10}{12.0}{\familydefault}{\mddefault}{\updefault}{\color[rgb]{0,0,0}$s_1=0$}%
}}}}
\put(4981,-772){\makebox(0,0)[lb]{\smash{{\SetFigFont{7}{8.4}{\familydefault}{\mddefault}{\updefault}{\color[rgb]{0,0,0}$(0,0,\kappa_1,0,\kappa_2-\kappa_1,0)$}%
}}}}
\put(2225,-4458){\makebox(0,0)[lb]{\smash{{\SetFigFont{10}{12.0}{\familydefault}{\mddefault}{\updefault}{\color[rgb]{0,0,0}$q=0$}%
}}}}
\put(1856,-6769){\makebox(0,0)[lb]{\smash{{\SetFigFont{12}{14.4}{\familydefault}{\mddefault}{\updefault}{\color[rgb]{0,0,0}$\mathbb{F}\cap(x_0=0)$}%
}}}}
\put(7046,-6769){\makebox(0,0)[lb]{\smash{{\SetFigFont{12}{14.4}{\familydefault}{\mddefault}{\updefault}{\color[rgb]{0,0,0}$\mathbb{F}$}%
}}}}
\put(5665,-38){\makebox(0,0)[lb]{\smash{{\SetFigFont{10}{12.0}{\familydefault}{\mddefault}{\updefault}{\color[rgb]{0,0,0}$s_0=s_1=0$}%
}}}}
\put(5468,-4998){\makebox(0,0)[lb]{\smash{{\SetFigFont{7}{8.4}{\familydefault}{\mddefault}{\updefault}{\color[rgb]{0,0,0}$(0,\kappa_1,0,0,\kappa_2,0)$}%
}}}}
\put(7525,-3876){\makebox(0,0)[lb]{\smash{{\SetFigFont{7}{8.4}{\familydefault}{\mddefault}{\updefault}{\color[rgb]{0,0,0}$(0,\kappa_1,0,\kappa_2,0,0)$}%
}}}}
\put(7805,-4998){\makebox(0,0)[lb]{\smash{{\SetFigFont{7}{8.4}{\familydefault}{\mddefault}{\updefault}{\color[rgb]{0,0,0}$(0,\kappa_1,0,0,0,\kappa_2)$}%
}}}}
\put(7319,-5926){\makebox(0,0)[lb]{\smash{{\SetFigFont{10}{12.0}{\familydefault}{\mddefault}{\updefault}{\color[rgb]{0,0,0}$s_0=q=0$}%
}}}}
\put(  1,-2308){\makebox(0,0)[lb]{\smash{{\SetFigFont{7}{8.4}{\familydefault}{\mddefault}{\updefault}{\color[rgb]{0,0,0}$(\kappa_1,0,0,0,\kappa_2,0)$}%
}}}}
\put(1403,-1653){\makebox(0,0)[lb]{\smash{{\SetFigFont{7}{8.4}{\familydefault}{\mddefault}{\updefault}{\color[rgb]{0,0,0}$(0,0,\kappa_1,0,\kappa_2-\kappa_1,0)$}%
}}}}
\put(3554,-1653){\makebox(0,0)[lb]{\smash{{\SetFigFont{7}{8.4}{\familydefault}{\mddefault}{\updefault}{\color[rgb]{0,0,0}$(0,0,\kappa_1,0,0,\kappa_2-\kappa_1)$}%
}}}}
\put(2439,-2446){\makebox(0,0)[lb]{\smash{{\SetFigFont{7}{8.4}{\familydefault}{\mddefault}{\updefault}{\color[rgb]{0,0,0}$(\kappa_1,0,0,0,0,\kappa_2)$}%
}}}}
\put(842,-3523){\makebox(0,0)[lb]{\smash{{\SetFigFont{7}{8.4}{\familydefault}{\mddefault}{\updefault}{\color[rgb]{0,0,0}$(0,\kappa_1,0,0,\kappa_2,0)$}%
}}}}
\put(3647,-3523){\makebox(0,0)[lb]{\smash{{\SetFigFont{7}{8.4}{\familydefault}{\mddefault}{\updefault}{\color[rgb]{0,0,0}$(0,\kappa_1,0,0,0,\kappa_2)$}%
}}}}
\end{picture}%

%% file: fibers2.pstex_t
\begin{picture}(0,0)%
\includegraphics{fibers2.pstex}%
\end{picture}%
\setlength{\unitlength}{2842sp}%
\begingroup\makeatletter\ifx\SetFigFont\undefined%
\gdef\SetFigFont#1#2#3#4#5{%
  \reset@font\fontsize{#1}{#2pt}%
  \fontfamily{#3}\fontseries{#4}\fontshape{#5}%
  \selectfont}%
\fi\endgroup%
\begin{picture}(8659,7276)(122,-7694)
\put(7876,-5161){\makebox(0,0)[lb]{\smash{{\SetFigFont{12}{14.4}{\familydefault}{\mddefault}{\updefault}{\color[rgb]{0,0,0}$X_{\infty}$}%
}}}}
\put(1876,-961){\makebox(0,0)[lb]{\smash{{\SetFigFont{7}{8.4}{\familydefault}{\mddefault}{\updefault}{\color[rgb]{0,0,0}$v=(0,0,\kappa_1,0,0,\kappa_2-\kappa_1)$}%
}}}}
\put(181,-3121){\makebox(0,0)[lb]{\smash{{\SetFigFont{9}{10.8}{\familydefault}{\mddefault}{\updefault}{\color[rgb]{0,0,0}$X_0\cap(q=0)=\mathbb{S}$}%
}}}}
\put(1981,-4621){\makebox(0,0)[lb]{\smash{{\SetFigFont{9}{10.8}{\familydefault}{\mddefault}{\updefault}{\color[rgb]{0,0,0}$X_{t_1}\cap(q=0)=\mathbb{S}$}%
}}}}
\put(4351,-2461){\makebox(0,0)[lb]{\smash{{\SetFigFont{7}{8.4}{\familydefault}{\mddefault}{\updefault}{\color[rgb]{0,0,0}$v$}%
}}}}
\put(4081,-6121){\makebox(0,0)[lb]{\smash{{\SetFigFont{9}{10.8}{\familydefault}{\mddefault}{\updefault}{\color[rgb]{0,0,0}$X_{t_2}\cap(q=0)=\mathbb{S}$}%
}}}}
\put(6451,-3961){\makebox(0,0)[lb]{\smash{{\SetFigFont{7}{8.4}{\familydefault}{\mddefault}{\updefault}{\color[rgb]{0,0,0}$v$}%
}}}}
\put(8551,-5461){\makebox(0,0)[lb]{\smash{{\SetFigFont{7}{8.4}{\familydefault}{\mddefault}{\updefault}{\color[rgb]{0,0,0}$v$}%
}}}}
\put(6526,-7636){\makebox(0,0)[lb]{\smash{{\SetFigFont{9}{10.8}{\familydefault}{\mddefault}{\updefault}{\color[rgb]{0,0,0}$\mathbb{S}$}%
}}}}
\put(1876,-586){\makebox(0,0)[lb]{\smash{{\SetFigFont{12}{14.4}{\familydefault}{\mddefault}{\updefault}{\color[rgb]{0,0,0}$X_0$}%
}}}}
\put(3676,-2161){\makebox(0,0)[lb]{\smash{{\SetFigFont{12}{14.4}{\familydefault}{\mddefault}{\updefault}{\color[rgb]{0,0,0}$X_{t_1}$}%
}}}}
\put(5776,-3661){\makebox(0,0)[lb]{\smash{{\SetFigFont{12}{14.4}{\familydefault}{\mddefault}{\updefault}{\color[rgb]{0,0,0}$X_{t_2}$}%
}}}}
\end{picture}%

%% file: P1xP2.pstex_t
\begin{picture}(0,0)%
\includegraphics{P1xP2.pstex}%
\end{picture}%
\setlength{\unitlength}{2842sp}%
\begingroup\makeatletter\ifx\SetFigFont\undefined%
\gdef\SetFigFont#1#2#3#4#5{%
  \reset@font\fontsize{#1}{#2pt}%
  \fontfamily{#3}\fontseries{#4}\fontshape{#5}%
  \selectfont}%
\fi\endgroup%
\begin{picture}(4907,4395)(976,-4690)
\put(1801,-3361){\makebox(0,0)[lb]{\smash{{\SetFigFont{11}{13.2}{\familydefault}{\mddefault}{\updefault}{\color[rgb]{0,0,0}$\mu$}%
}}}}
\put(2401,-511){\makebox(0,0)[lb]{\smash{{\SetFigFont{12}{14.4}{\familydefault}{\mddefault}{\updefault}{\color[rgb]{0,0,0}$S_{\infty}$}%
}}}}
\put(976,-736){\makebox(0,0)[lb]{\smash{{\SetFigFont{12}{14.4}{\familydefault}{\mddefault}{\updefault}{\color[rgb]{0,0,0}$S_{\infty}'$}%
}}}}
\put(3526,-4636){\makebox(0,0)[lb]{\smash{{\SetFigFont{11}{13.2}{\familydefault}{\mddefault}{\updefault}{\color[rgb]{0,0,0}$1$}%
}}}}
\put(5101,-2986){\makebox(0,0)[lb]{\smash{{\SetFigFont{12}{14.4}{\familydefault}{\mddefault}{\updefault}{\color[rgb]{0,0,0}$L^{\mu}$}%
}}}}
\put(3076,-1186){\makebox(0,0)[lb]{\smash{{\SetFigFont{10}{12.0}{\familydefault}{\mddefault}{\updefault}{\color[rgb]{0,0,0}$x_2=0$}%
}}}}
\put(5176,-3811){\makebox(0,0)[lb]{\smash{{\SetFigFont{12}{14.4}{\familydefault}{\mddefault}{\updefault}{\color[rgb]{0,0,0}$S_0$}%
}}}}
\put(1276,-4261){\makebox(0,0)[lb]{\smash{{\SetFigFont{10}{12.0}{\familydefault}{\mddefault}{\updefault}{\color[rgb]{0,0,0}$s_0=x_1=0$}%
}}}}
\put(4651,-2161){\makebox(0,0)[lb]{\smash{{\SetFigFont{10}{12.0}{\familydefault}{\mddefault}{\updefault}{\color[rgb]{0,0,0}$s_1=x_0=0$}%
}}}}
\end{picture}%